\documentclass[twocolumn  
]{svjour3} 

\usepackage{graphicx,times}

\usepackage[breaklinks=true, colorlinks=true, pdfstartview=FitV, linkcolor=blue, citecolor=blue, urlcolor=red]{hyperref}

\usepackage{amsmath}
\usepackage{amsfonts}

\allowdisplaybreaks

\journalname{Nonlinear Dynamics}

\begin{document}

\title{Solution to the main problem of the artificial satellite by reverse normalization}

\author{Martin Lara}

\institute{M. Lara \at
              Scientific Computing Group--GRUCACI, University of La Rioja, Madre de Dios 53, 26006 Logro\~no, La Rioja, Spain   \\
              Tel.: +34-941-299440 \\
              Fax: +34-941-299460 \\
              \email{mlara0@gmail.com}           
           \and
           M. Lara \at
              Space Dynamics Group, Polytechnic University of Madrid--UPM, Plaza Cardenal Cisneros 3, 28040 Madrid, Spain
}

\date{\color{red} draft of April 16, 2020}
\maketitle

\begin{abstract}
The non-linearities of the dynamics of Earth artificial satellites are encapsulated by two formal integrals that are customarily computed by perturbation methods. Standard procedures begin with a Hamiltonian simplification that removes non-essential short-period terms from the Geopotential, and follow with the removal of both short- and long-period terms by means of two different canonical transformations that can be carried out in either order. We depart from the tradition and proceed by standard normalization to show that the Hamiltonian simplification part is dispensable. Decoupling first the motion of the orbital plane from the in-plane motion reveals as a feasible strategy to reach higher orders of the perturbation solution, which, besides, permits an efficient evaluation of the long series that comprise the analytical solution.

\keywords{main problem \and Hamiltonian mechanics \and normalization \and canonical perturbation theory \and floating-point arithmetic}
\end{abstract}

\section{Introduction}

The dynamics of close Earth satellites under gravitational effects are mostly driven by perturbations of the Keplerian motion induced by the Earth oblateness. For this reason, the approximation obtained when truncating the Legendre polynomials expansion of the Geopotential to the only contribution of the zonal harmonic of the second degree, whose coefficient is usually denoted $J_2$, is traditionally known as the \emph{main problem} of artificial satellite theory. The $J_2$-truncation of the gravitational potential is known to give rise to nonintegrable dynamics \cite{Danby1968,IrigoyenSimo1993,Broucke1994,CellettiNegrini1995} that comprise short- and long-period effects, as well as secular terms \cite{Kozai1959,Kaula1961}. However, due to the smallness of the $J_2$ coefficient of the Earth, the full system can be replaced by a separable approximation, which is customarily obtained by removing the periodic effects by means of perturbation methods \cite{Nayfeh2004}.
\par

When written in the action-angle variables of the Kepler problem, also called Delaunay variables, the main problem Hamiltonian immediately shows that the right ascension of the ascending node is a cyclic variable. In consequence, its conjugate momentum, the projection of the angular momentum vector along the Earth's rotation axis, is an integral of the main problem, which, therefore, is just of two degrees of freedom. Then, following Brouwer \cite{Brouwer1959}, the main problem Hamiltonian is normalized in two steps. First, the short period effects are removed by means of a canonical transformation that, after truncation to some order of the perturbation approach, turns the conjugate momentum to the mean anomaly (the Delaunay action) into a formal integral. Next, a new canonical transformation converts the total angular momentum into another formal integral. The main problem Hamiltonian is in this way completely reduced to a function of only the momenta in the new variables, whose Hamilton equations are trivially integrable.
\par

In spite of the normalization procedure is standard from the point of view of perturbation theory, it happens that not all the action-angle variables of the Kepler problem appear explicitly in the Geopotential, as is usually advisable in the construction of perturbation solutions \cite{Kinoshita1972,LaraFerrer2013,LaraFerrer2015,Lara2018}. In particular, the mean anomaly remains implicit in the gravitational potential through its dependence on the polar angle. Unavoidably, Kepler's equation must be solved to show the mean anomaly explicit, a fact that entails expanding  the elliptic motion in powers of the eccentricity \cite{Hansen1855,BrouwerClemence1961}. This standard procedure is quite successful when dealing with orbits of low eccentricities, like is typical in a variety of astronomy problems \cite{Delaunay1860,Tisserand1889}. However, it involves handling very long Fourier series in the case of orbits with moderate eccentricities \cite{DepritRom1970,Kinoshita1977SAOSR}, and hence the application of this method to different problems of interest in astrodynamics is de facto prevented.
\par

On the contrary, the integrals appearing in the solution of the artificial satellite problem can be solved in closed form of the eccentricity \cite{Brouwer1959,Kozai1959}. It only requires the help of the standard relation between the differentials of the true and mean anomalies that is derived from the preservation of the total angular momentum of the Keplerian motion. Regrettably, the closed form approach soon finds difficulties in achieving higher orders of the short-period elimination, which stem from the impossibility of obtaining the antiderivative of the equation of the center (the difference between the true and mean anomalies) in closed form of the eccentricity in the realm of trigonometric functions \cite{Jefferys1971,OsacarPalacian1994}. Nonetheless, the difficulties are overcome by the artifact of grouping the equation of the center with other functions appearing in the procedure previously to approaching their integration \cite{Kozai1962,Aksnes1971,Deprit1982}. Alternatively, the application of a preliminary Hamiltonian simplification, the \emph{elimination of the parallax} \cite{Deprit1981,LaraSanJuanLopezOchoa2013b}, eases the consequent removal of short-period effects to some extent \cite{CoffeyDeprit1982,Healy2000}.
\par

From the point of view of the perturbation approach, removing the short-period effects in the first place seems the more natural in view of the degeneracy of the Kepler problem. Indeed, the Kepler Hamiltonian in action-angle variables, on which the perturbation approach hinges on, only depends on the Delaunay action \cite{GoldsteinPooleSafko2001,LaraTossa2016}. However, the order in which the formal integrals are sequentially introduced in the perturbation solution is not relevant in a total normalization procedure, whose result is unique \cite{Arnold1989}. In fact, it happens that relegating the transformation of the Delaunay action into a formal integral to the last step of the perturbation approach provides clear simplifications in dealing with the equation of the center \cite{AlfriendCoffey1984}. In this way the task of extending the solution of the main problem to higher orders is notably simplified.
\par

Converting the total angular momentum into a formal integral of the main problem requires making cyclic the argument of the perigee, up to some truncation order of the perturbation approach, in the transformed Hamiltonian. However, as it appeared in the literature, the transformation called by their authors the \emph{elimination of the perigee} \cite{AlfriendCoffey1984} is not the typical normalization procedure, although it operates analogous results. Indeed, on the one hand, the elimination of the perigee is only applied to a Hamiltonian obtained after the elimination of the parallax, to which simplification it could seem to be unavoidably attached. On the other hand, when extended to higher orders, it removes more terms than those needed in a partial normalization. Last, the fact that the technique was originally devised in the canonical set of polar variables, to which the argument of the perigee does not pertain, neither helps in grasping the essence of the transformation. Reimplementation of the procedure in the usual set of action-angle variables makes the process of converting the argument of the perigee into a cyclic variable much more evident \cite{LaraSanJuanLopezOchoa2013c}, but it still bears the same differences with respect to a classical normalization procedure.
\par

We disregard the claimed benefits of Hamiltonian simplification procedures and compute the solution to the main problem of the artificial satellite theory by standard normalization. It is called \emph{reverse} normalization because we exchange the order in which the formal integrals are traditionally introduced when solving the artificial satellite problem. More precisely, the total angular momentum is transformed into a formal integral in the first place, in this way decoupling the motion of the orbital plane from the satellite's motion on that plane.\footnote{The advantages of decoupling the motion of the instantaneous orbital plane from the in-plane motion are well known, and are commonly pursued in the search for efficient numerical integration methods, q.v.~\cite{Lara2017if} and references therein.} Then, a second canonical transformation converts the mean anomaly into a cyclic variable, in this way achieving the total reduction of the main problem Hamiltonian. 
\par

The procedure for making the argument of the perigee cyclic in the first place follows an analogous strategy to the one devised in the classical elimination of the perigee transformation \cite{AlfriendCoffey1984}. However, in our approach it is applied directly to the original main problem Hamiltonian, and differs from the original technique, as well as from an analogous procedure carried out in \cite{SanJuanetal2013}, in which the parallactic terms (inverse powers of the radius with exponents higher than 2) are not removed from the new, partially normalized Hamiltonian. In spite of that, we did not find trouble in dealing with the equation of the center in closed form in the subsequent Delaunay normalization \cite{Deprit1982}, a convenience that might had been anticipated from the discussions in \cite{Lara2019CMDA}.
\par

The Hamiltonian reduction has been approached in Delaunay variables. Unfortunately, these variables share the deficiencies of their partner Keplerian elements, which are singular for circular orbits and for equatorial orbits. Because of that, the secular terms are reformulated in a set of non-singular variables that replaces the mean anomaly, the argument of the perigee, and the total angular momentum, by the mean argument of the latitude and the projections of the eccentricity vector in the nodal frame, which are sometimes denoted \emph{semi-equinoctial} variables \cite{Konopliv1990}. For the periodic corrections, we find convenience in using polar variables, which in the particular case of the main problem are free from small divisors of any kind except for those related to the critical inclination resonance, and are known to produce compact expressions that yield faster evaluation \cite{Izsak1963,LaraVilhenaSanchezPrado2015}.
\par

We extended the complete normalization to the order six of the perturbation approach, which, to our knowledge, is the maximum order that has been reported in the literature (yet limited to partial normalization cases) \cite{Healy2000,SanJuanetal2013}. The aim of computing such a high order is not to enter a competition. On the contrary, we did it simply because the particular value of the Earth's $J_2$ coefficient is $\sim10^{-3}$, and hence the computed solution should be exact to the numerical precision of standard floating-point arithmetic \cite{Kahan1997} already at the fifth order. We checked that it is exactly the case, and extending the computations to the sixth order only served us to verify that we don't find observable improvements in our tests. In order to compute this high-order solution we resorted to the practicalities of standard commercial software, in which Deprit's perturbation algorithm based on the Lie transforms method \cite{Deprit1969} is easily implemented. On account of the current computational power, it should be quite feasible to extend the perturbation solution, if desired, to even higher orders, although we didn't try that.
\par

On the other hand, making cyclic the argument of the perigee with the standard normalization in action-angle variables seems to be a particularly efficient procedure from the computational point of view, resulting in a generating function whose size is astonishingly smaller than alternative proposals in the literature. This fact, combined with the reduction in the number of transformations required by the perturbation solution to just two, as opposite to the three transformations needed when the elimination of the parallax is carried out in the first place, might make this latter simplification dispensable, as well as the discussions in \cite{DepritMiller1989} questionable, at least in what respects to the computation of explicit analytical solutions as opposite to partial normalizations to be integrated semianalytically. 
\par

Tests carried out on typical Earth orbits of different types confirm that the computed solution bears exactly the expected characteristics of a perturbation solution of perturbed Keplerian motion. In particular, we verified the degree of convergence of successive approximations. We confirmed too that, as expected, the quality of the analytical solution degrades in the vicinity of the critical inclination resonance ---a case that, of course, would require a specific treatment \cite{CoffeyDepritMiller1986,Jupp1988,Lara2015IR}. In addition, we checked that the computation of the initialization constants of the analytical solution from corresponding initial conditions gets a clear benefit of using a higher order of the periodic corrections than the order needed for ephemeris evaluation, yet this additional accuracy is not needed in all the variables and can be limited to the computation of the formal integral given by the Delaunay action \cite{BreakwellVagners1970,DepritRom1970,HautesserresLaraCeMDA2017}.

\section{Construction of the analytical solution}

The solution of Laplace's equation in spherical coordinates provides the gravitational potential in the form of an expansion in harmonic functions. In the case of the Earth, the mass distribution is almost symmetric with respect to the rotation axis. For this reason, the Geopotential is simplified in different applications to just the contribution of the zonal harmonics (see \cite{Danby1992}, for instance) 
\begin{equation} \label{zonalpot}
\mathcal{V}=\frac{\mu}{r}\sum_{m\ge{0}}\frac{R_\oplus^m}{r^m}J_mP_{m}(\sin\varphi),
\end{equation}
where $r$ is distance from the origin, $\varphi$ is latitude, $\mu$ is the Earth's gravitational parameter, $R_\oplus$ is the Earth's mean equatorial radius, $J_m$ are zonal harmonic coefficients, and $P_{m}$ denote Legendre polynomials of degrees $m$.
\par

The flow stemming from the potential (\ref{zonalpot}) admits Hamiltonian formulation. Besides, because $J_2\sim10^{-3}$ whereas $J_m$ ($m>2$) are $\mathcal{O}(J_2^2)$ or smaller, the Hamiltonian
\begin{equation} \label{Hammain}
\mathcal{H}=-\frac{\mu}{2a}-\frac{\mu}{r}\frac{R_\oplus^2}{r^2}J_2\frac{1}{2}\left(1-3\sin^2\varphi\right),
\end{equation}
is representative of the main characteristics of the dynamics of close Earth satellites, and is commonly known as the \emph{main problem} of artificial satellite theory. The symbol $a$ in Eq.~(\ref{Hammain}) stands for the orbit semi-major axis, and, from standard trigonometric relations, $\sin\varphi=\sin{I}\sin\theta$, where $I$ and $\theta$ are used to denote the inclination and the argument of the latitude, respectively.
\par

When using the classic set of Keplerian elements, the argument of the latitude is $\theta=f+\omega$, where $f$ and $\omega$ denote the true anomaly and the argument of the periapsis, respectively. The radius is computed from the conic equation
\begin{equation} \label{conic}
r=\frac{p}{1+e\cos{f}},
\end{equation}
where $p=a(1-e^2)$ is the parameter of the conic and $e$ is the eccentricity of the orbit. Note, however, that, because we are using Hamiltonian formulation, the symbols appearing in Eq.~(\ref{Hammain}) need to be expressed as functions of some set of canonical variables. While the Keplerian variables are not canonical, they are conveniently expressed in the set of Delaunay canonical variables $(\ell,g,h,L,G,H)$, in which $\ell$ is the mean anomaly, $g=\omega$, $h$ is the argument of the ascending node, $L=\sqrt{\mu{a}}$ is known as the Delaunay action, $G=L\sqrt{1-e^2}$ is the total angular momentum, and $H=G\cos{I}$ is the projection of the angular momentum along the Earth's rotation axis. That the latter is an integral of the main problem becomes evident after checking that $h$ is an ignorable variable in Hamiltonian (\ref{Hammain}), which, in consequence, is of just two degrees of freedom.
\par

Besides, the main problem Hamiltonian is conservative, and, therefore, Eq.~(\ref{Hammain}) remains constant (the total energy) for a given set of initial conditions. On the other hand, the existence of a third integral has not been proved. Then, exact solutions of the main problem are not known. Alternatively, it is approached with the usual tools of non-linear dynamics, such as Poincar\'e surfaces of section or the computation of families of periodic orbits \cite{Danby1968,Broucke1994,SADSaM1999,Lara1999}. Still, in some regions of phase space, and for particular energy values, the third integral can be constructed formally with the help of perturbation methods, obtaining in this way useful analytical approximations to the main problem dynamics in these particular regions.
\par

\subsection{Perturbation approach}

We rewrite the main problem in the form of a perturbation Hamiltonian
\begin{equation} \label{mainP}
\mathcal{H}=\mathcal{H}(\ell,g,-,L,G,H)\equiv\sum_{m\ge0}\frac{\varepsilon^m}{m!}K_{m,0},
\end{equation}
in which $\varepsilon$ is a formal small parameter ($\varepsilon=1$) used to manifest the strength of each summand of Eq.~(\ref{mainP}), and
\begin{eqnarray} \label{mainPK0}
K_{0,0} &=& -\frac{\mu}{2a},  \\ \label{mainPK1}
K_{1,0} &=& -\frac{\mu}{r}\frac{R_\oplus^2}{r^2}\frac{1}{4}J_2\left[2-3s^2+3s^2\cos(2f+2\omega)\right], \\ \label{mainPKm}
K_{m,0} &=& 0, \quad m\ge2.
\end{eqnarray}
The symbol $s$ in Eq.~(\ref{mainPK1}) stands for the usual abbreviation of the sine of the inclination.
\par

The aim is to find a transformation of variables
\[
(\ell,g,h,L,G,H;\varepsilon)\mapsto(\ell',g',h',L',G',H'),
\]
given also by an expansion in powers of the small parameter, such that, up to some truncation order $m=k$, the Hamiltonian in the new variables is transformed into a separable Hamiltonian. Namely,
\begin{equation} \label{mainPp}
\mathcal{H}'=\sum_{m=0}^{k}\frac{\varepsilon^m}{m!}K_{0,m}(-,-,-,L',G',H)+\mathcal{O}(\varepsilon^{k+1}).
\end{equation}
\par

The desired transformation is derived from a scalar generating function 
\begin{equation} \label{Lie:Gen}
\mathcal{W}=\sum_{m\ge0}\frac{\varepsilon^m}{m!}W_{m+1}.
\end{equation}
which, in our case, is computed using Deprit's perturbation algorithm based on Lie transforms \cite{Deprit1969}. The procedure is summarized in finding a particular solution of the \emph{homological equation}
\begin{equation} \label{homoeq}
\mathcal{L}_0(W_m)\equiv\{W_m;K_{0,0}\}=\widetilde{K}_{0,m}-K_{0,m},
\end{equation}
in which $\{W_m;K_{0,0}\}$ stands for the computation of the Poisson bracket of $W_m$ and the zeroth order term of the Hamiltonian. Terms
$\widetilde{K}_{0,m}$ are known from the original Hamiltonian as well as from previous computations to the step $m$, which are carried out using Deprit's fundamental recursion
\begin{equation} \label{Lie:triangle}
{F}_{n,q}={F}_{n+1,q-1}+
\sum_{m=0}^n{n\choose{m}}\{{F}_{n-m,q-1};W_{m+1}\}.
\end{equation}
Equation (\ref{Lie:triangle}) can be used to formulate any function $\mathcal{F}$ of the canonical set of original variables like a function of the new variables, the Hamiltonian (\ref{mainP}) being just one among the different possibilities. Finally, the term $K_{0,m}$ is chosen at our will, with the only condition of making the homological equation solvable in $W_m$.
\par

When using Delaunay variables the homological equation of the main problem Hamiltonian is solved by indefinite integration. Indeed, Eq.~(\ref{mainPK0}) turns into $K_{0,0}=-\frac{1}{2}\mu^2/L^2$, and hence the \emph{Lie derivative} operator $\mathcal{L}_0$ takes the extremely simple form
\[
\mathcal{L}_0(W_m)=n\frac{\partial{W}_m}{\partial\ell},
\]
in which $n=\sqrt{\mu/a^3}=\mu^2/L^3$ is the mean motion of the unperturbed problem. Then, from Eq.~(\ref{homoeq}),
\begin{equation} \label{homosoll}
W_m=\frac{1}{n}\int(\widetilde{K}_{0,m}-K_{0,m})\,\mathrm{d}\ell+C_m,
\end{equation}
where the functions $C_m\equiv{C}_m(-,g,L,G;H)$ play the role of integration ``constants'' that verify $\mathrm{d}C_m/\mathrm{d}\ell=0$. Equation (\ref{homosoll}) is solved in closed form of the eccentricity with the help of the differential relation between the true and mean anomalies provided by the preservation of the angular momentum of the Kepler problem. That is, $G=r^2\mathrm{d}f/\mathrm{d}t$, from which
\begin{equation} \label{dl2df}
a^2\eta\,\mathrm{d}\ell=r^2\mathrm{d}f,
\end{equation}
where $\eta=\sqrt{1-e^2}=G/L$ is usually known as the \emph{eccentricity function}. Hence,\begin{equation} \label{homosol}
W_m=\frac{1}{n}\int(\widetilde{K}_{0,m}-K_{0,m})\frac{r^2}{a^2\eta}\,\mathrm{d}f+C_m.
\end{equation}
\par

The transformation that makes the main problem Hamiltonian separable, up to the truncation order, is split into two different canonical transformations. With the first one 
\[
(\ell,g,h,L,G,H;\varepsilon)\mapsto(\ell',g',h',L',G',H'),
\]
we make the argument of the perigee cyclic, thus converting the total angular momentum into a formal integral. That is, the motion of the satellite in the orbital plane, whose inclination remains constant in the new variables, is decoupled from the motion of the orbital plane. Therefore, the reduced system representing the motion on the orbital plane could be integrated separately. Rather, we carry out a second canonical transformation
\[
(\ell',g',h',L',G',H';\varepsilon)\mapsto(\ell'',g'',h'',L'',G'',H''),
\]
that makes ignorable the mean anomaly in the transformed Hamiltonian in double-prime variables. In this way, the Delaunay action is converted into a formal integral too, and the complete Hamiltonian reduction is achieved up to the truncation order of the perturbation solution, whose Hamilton equations are thus trivially integrable.
\par

The computation of the secular terms from the normalized Hamiltonian is only part of the perturbation solution. To be complete, it requires the correct initialization of the constants of the perturbation solution from corresponding initial conditions using the \emph{inverse} transformation (from original to double prime variables), on the one hand, and the recovery of periodic effects using the \emph{direct} transformation (from double prime to original variables) to obtain the ephemeris corresponding to the secular terms propagation. Both the direct and inverse transformations are obtained from standard application of the Lie transforms method. Because there are no specific artifices related to the computation of the current solution in what respects to that part, we do not provide details on their computation and an interested reader is referred to the literature \cite{Deprit1969}.

\subsection{Normalization of the total angular momentum}

At the first order, we check from Eq.~(\ref{Lie:triangle}) that the known terms involved in the homological equation consist only of $\widetilde{K}_{0,1}=K_{1,0}$. Therefore, we chose the new Hamiltonian term 
\begin{equation} \label{K01g}
K_{0,1}=-\frac{\mu}{r}\frac{R_\oplus^2}{r^2}\frac{1}{4}J_2\left(2-3s^2\right),
\end{equation}
which is the part of Eq.~(\ref{mainPK1}) that is free from the argument of the perigee, as desired. Then, Eq.~(\ref{homosol}) turns into
\[
W_1=-\frac{3}{4}GJ_2\frac{R_\oplus^2}{p^2}s^2\int\frac{p}{r}\cos(2f+2\omega)\,\mathrm{d}f+C_1,
\]
which is solved by standard integration after replacing the radius with the conic equation (\ref{conic}). We obtain
\begin{equation} \label{W1g}
W_1=-\epsilon G \frac{1}{2}s^2\sum_{i=1}^3 3^{\left\lfloor2-\frac{1}{2}i\right\rfloor } e^{\left| i-2\right| } \sin (if+2\omega)+C_1,
\end{equation}
where $C_1$ is hold arbitrary by now, the symbol $\lfloor\;\rfloor$ stands for integer part, and we abbreviated 
\begin{equation} \label{mieps}
\epsilon=\frac{1}{4}J_2\frac{R_\oplus^2}{p^2},
\end{equation}
in which $p=G^2/\mu$, and, therefore $\epsilon$ is a function of the total angular momentum. Recall that the symbols $p$, $s$, $e$, $\omega$, and $f$ in Eqs.~(\ref{K01g}), (\ref{W1g}), and (\ref{mieps}) are functions of the Delaunay variables.
\par

At the second order we compute $K_{0,2}$ from Eq.~(\ref{Lie:triangle}), to obtain
\begin{equation} \label{triangle2}
K_{0,2}=\{K_{0,1},W_1\}+K_{1,1},
\end{equation}
in which $K_{1,1}$ is computed using again Eq.~(\ref{Lie:triangle}). Namely,
\begin{equation} \label{K11}
K_{1,1}=\{K_{0,0},W_2\}+\{K_{1,0},W_1\}+K_{2,0}.
\end{equation}
On account of $K_{2,0}=0$ from Eq.~(\ref{mainPKm}), the known terms hitherto of the homological equation (\ref{homoeq}) are
\[
\widetilde{K}_{0,2}=\{K_{0,1},W_1\}+\{K_{1,0},W_1\},
\]
whose computation only involves partial differentiation. Before approaching the solution of Eq.~(\ref{homosol}), we first check that the term $\widetilde{K}_{0,2}$ is made of trig\-o\-nom\-e\-tric coefficients $T_i$ whose arguments always involve the true anomaly as argument, except for the terms
\begin{eqnarray*}
T_0 &=& \epsilon^2\frac{3}{4}\frac{\mu}{r}\frac{p}{r}s^2\left[e^2(23s^2-16\right)+8 \left(4 s^2-3)\right],
\\
T_1 &=&-\frac{3}{2}\epsilon^2G^2\frac{1}{r^2}(15s^2-14)s^2e^2\cos2\omega,
\\
T_2 &=& 6\epsilon{G}\frac{1}{r^2}(5s^2-4)\frac{\partial{C}_1}{\partial{g}}.
\end{eqnarray*}
Integration of these 3 terms in Eq.~(\ref{homosol}) would yield corresponding terms that grow unbounded with $f$.
The term $T_0$ is free from the argument of the periapsis, and hence it can be cancelled by choosing the new Hamiltonian $K_{0,2}$ having a corresponding term $T_0$. On the contrary, the terms $T_1$ and $T_2$ depend on the argument of the perigee, a fact that prevents their appearance in the new Hamiltonian under our requirement of making $g$ a cyclic variable. Nevertheless, since we had left $C_1$ arbitrary, both terms will cancel each other if $C_1$ is determined from the partial differential equation $T_1+T_2=0$. We readily check that the particular solution
\begin{equation} \label{C1g}
C_1=\epsilon{G}\frac{15s^2-14}{8(5s^2-4)}s^2e^2\sin2\omega,
\end{equation}
matches this purpose. Remark that the appearance of the divisor $5s^2-4$ prevents application of the solution to the resonant cases that happen at the so-called \emph{critical inclinations} in which $\sin^2I=4/5$. That is, the supplementary inclinations $I=63.435^\circ$ and $I=116.565^\circ$.
\par

After computing the partial derivatives of Eq.~(\ref{C1g}) with respect to $L$ and $G$, which also appear among the coefficients of the terms remaining in $\widetilde{K}_{0,2}$, the known terms of the homological equation become fully determined. Then, we can safely choose $K_{0,2}$ so that it comprises those terms of $\widetilde{K}_{0,2}$ that are free from the argument of the periapsis. We obtain,
\begin{equation} \label{K02g}
K_{0,2}=\epsilon^2\frac{\mu}{r}\frac{p^2}{r^2}\frac{3s^2}{8(5s^2-4)^2}\sum_{j=0}^2\frac{p^j}{r^j}\sum_{k=0}^{\lfloor1-\frac{1}{2}j\rfloor}e^{2k}\gamma_{2,j,k},
\end{equation}
in which the inclination polynomials $\gamma_{2,j,k}\equiv\gamma_{2,j,k}(s)$ are provided in Table~\ref{t:K2gpoly}. 
\par

\begin{table}[htb]
\caption{Inclination polynomials $\gamma_{2,j,k}$ in Eq.~(\protect\ref{K02g}). }
\label{t:K2gpoly}
\begin{tabular}{@{}ll@{}}
\hline\noalign{\smallskip}
${}_{0,0}$ & $-8 \left(200 s^6-455 s^4+345 s^2-88\right)$ \\
${}_{0,1}$ & $375 s^6-930 s^4+780s^2-224$ \\
${}_{1,0}$ & $5 \left(805 s^6-1878 s^4+1464 s^2-384\right)$ \\
${}_{2,0}$ & $-825 s^6+1990s^4-1616 s^2+448$ \\
\noalign{\smallskip}\hline
\end{tabular}
\end{table}

The second order term of the generating function is then readily computed from Eq.~(\ref{homosol}), yielding
\begin{eqnarray} \nonumber
W_2 &=& \frac{\epsilon^2G}{32(5s^2-4)^2}\sum_{l=1}^2\sum_{\substack{k=l-2\\k\ne0}}^{2l+2}\sum_{j=0}^1\Gamma_{2,j,k,l}e^{2j+k^*} \\ \label{W2g}
&& \times s^{2l}\sin(kf+2l\omega)+C_2,
\end{eqnarray}
with $k^*\equiv{k}\bmod2$, and the inclination polynomials $\Gamma_{2,j,k,l}$ are given in Table \ref{t:W2gpoly}.
\par

\begin{table}[htb]
\caption{Non-zero inclination polynomials $\Gamma_{2,j,k,l}$ in Eqs.~(\protect\ref{W2g})--(\protect\ref{C2g}). }
\label{t:W2gpoly}
\begin{tabular}{@{}ll@{}}
\hline\noalign{\smallskip}
${}_{1,-1,1}$ & $-12(5 s^2-4) (7 s^2-6) (15s^2-14)$ \\
${}_{0,1,1}$ & $-48(5 s^2-4) (195 s^4-340s^2+148)$ \\
${}_{1,1,1}$ & $24(5 s^2-4)^2 (15 s^2-14)$ \\
${}_{1,1,2}$ & $-3 (225 s^4-430 s^2+208)$ \\
${}_{0,2,1}$ & $-96 (5s^2-4)^2 (9 s^2-8)$ \\
${}_{1,2,1}$ & $-24 (5 s^2-4) (65s^4-116 s^2+52)$ \\
${}_{1,2,2}$ & $-60 (50 s^4-87 s^2+38)$ \\
${}_{0,3,1}$ & $-64 (5 s^2-4)^2 (8 s^2-7)$ \\
${}_{1,3,1}$ & $4 (3s^2-2) (5 s^2-4) (15 s^2-14)$ \\
${}_{0,3,2}$ & $-4 (5s^2-4) (135 s^2-122)$ \\
${}_{1,3,2}$ & $-8 (75 s^4-135s^2+61)$ \\
${}_{1,4,1}$ & $-12 (5 s^2-4)^2 (7 s^2-6)$ \\
${}_{0,4,2}$ & $24 (5 s^2-4)^2$ \\
${}_{1,4,2}$ & $-12 (5 s^2-4) (25s^2-23)$ \\
${}_{0,5,2}$ & $24 (5 s^2-4)^2$ \\
${}_{1,5,2}$ & $-3 (5s^2-4) (15 s^2-14)$ \\
${}_{1,6,2}$ & $6 (5 s^2-4)^2$ \\
${}_{0,0,2}$ & $(15 s^2-14)^2(15s^2-13)$ \\
${}_{0,0,1}$ & $8(5s^2-4)^2(1215s^4-1997s^2+824)$ \\
${}_{1,0,1}$ & $-2(5s^2-4)(15s^2-14)(45s^4+36s^2-56)$ \\
\noalign{\smallskip}\hline
\end{tabular}
\end{table}

Analogously as we did with $C_1$, the arbitrary function $C_2$ in which Eq.~(\ref{W2g}) depends upon is determined at the next step of the perturbation approach in such a way that the appearance of unbounded terms in the solution of the homological equation of the third order are avoided. Thus, we compute $\widetilde{K}_{0,3}$ by successive evaluations of the fundamental recursion (\ref{Lie:triangle}). After identifying the problematic terms of $\widetilde{K}_{0,3}$, they are cancelled by computing
\begin{equation} \label{C2g}
C_2=\frac{\epsilon^2G}{64(5s^2-4)^3}\sum_{l=1}^2\sum_{j=0}^{2-l}\Gamma_{2,j,0,l}e^{2(j+l)}s^{2l}\sin2l\omega,
\end{equation}
which contributes three additional inclination functions that are also listed in Table \ref{t:W2gpoly}. In this way, the second order term of the generating function given by Eq.~(\ref{W2g}) becomes fully determined.
\par

The needed partial derivatives of $C_2$ appearing in $\widetilde{K}_{0,3}$ are then computed, and the third order term of the new Hamiltonian is chosen, as before, to comprise those terms of $\widetilde{K}_{0,3}$ that are free from the argument of the periapsis. We obtain
\begin{equation} \label{K03g}
K_{0,3}=\epsilon^3\frac{\mu}{r}\frac{p^2}{r^2}\frac{3s^2}{32(5s^2-4)^3}\sum_{j=0}^4\frac{p^j}{r^j}\sum_{k=0}^{\lfloor2-\frac{1}{2}j\rfloor}e^{2k}\gamma_{3,j,k},
\end{equation}
where the inclination polynomials $\gamma_{3,j,k}$ are displayed in Table \ref{t:K3g}.
\par

\begin{table}[htb]
\caption{Inclination polynomials $\gamma_{3,j,k}$ in Eq.~(\protect\ref{K03g}). }
\label{t:K3g}
\begin{tabular}{@{}ll@{}}
\hline\noalign{\smallskip}
${}_{0,0}$ & $-8(5s^2-4)(313525 s^8-899030 s^6+933656 s^4$ \\
 & $-409296 s^2+61824)$ \\
${}_{0,1}$ & $4(1551625 s^{10}-5675700 s^8+8148960 s^6-5706408 s^4$ \\
 & $+1930272 s^2-248064)$ \\
${}_{0,2}$ & $-2 (40500 s^{10}-99525 s^8+64840 s^6+18788 s^4$ \\
 & $-33936 s^2+9408)$ \\
${}_{1,0}$ & $2 (5 s^2-4) (2631475 s^8-7558270 s^6+7872692 s^4$ \\
 & $-3470616 s^2+530304)$ \\
${}_{1,1}$ & $-3457125 s^{10}+12282750 s^8-17085020 s^6$ \\
 & $+11554040 s^4-3756000 s^2+459648$ \\
${}_{2,0}$ & $-2 (5 s^2-4) (1584375 s^8-4536150 s^6+4716436 s^4$ \\
 & $-2082712 s^2+321408)$ \\
${}_{2,1}$ & $138375 s^{10}-128250 s^8-351900 s^6+612440 s^4$ \\
 & $-326368 s^2+56448$ \\ 
${}_{3,0}$ & $8 (5 s^2-4) (93300 s^8-259915 s^6+264982 s^4-116928 s^2$ \\
 & $+18816)$ \\
${}_{4,0}$ & $-20 s^2 (5 s^2-4) (15 s^2-14) (45 s^4+36 s^2-56)$
\\
\noalign{\smallskip}\hline
\end{tabular}
\end{table}

Then, after solving Eq.~(\ref{homosol}), we find that the third order term of the generating function can be arranged in the form
\begin{eqnarray} \nonumber
W_3 &=& \frac{\epsilon^3G}{8960 (5 s^2-4)^4}\sum_{l=1}^3\sum_{\substack{k=l-4 \\ k\ne0}}^{2l+4}\sum_{j=0}^2\Gamma_{3,j,k,l}e^{2j+k^*} \\ \label{W3g}
&& \times s^{2l}\sin(kf+2lg)+C_3
\end{eqnarray}
in which the inclination polynomials $\Gamma_{3,j,k,l}$ that do not vanish are listed in Table \ref{t:W3g}. The constant $C_3$ is determined at the next order of the perturbation approach in the same way as we did previously. We obtain
\[
C_3=\frac{\epsilon^3G}{1536(5s^2-4)^5}\sum_{l=1}^3\sum_{j=0}^{3-l}\Gamma_{3,j,0,l}e^{2(j+l)}s^{2l}\sin2l\omega,
\]
which contributes six additional inclination polynomials $\Gamma_{3,j,0,l}$, that are also listed in Table \ref{t:W3g}.
\par

\begin{table*}[htbp]
\caption{Non-zero inclination polynomials $\Gamma_{3,j,k,l}$ in Eq.~(\protect\ref{W3g}). }
\label{t:W3g}
\begin{tabular}{@{}ll@{}}
\hline\noalign{\smallskip}
${}_{2,-3,1}$ & $35 \left(15 s^2-14\right) \left(87375 s^{10}-335550 s^8+505080
   s^6-371184 s^4+132096 s^2-17920\right)$ \\
${}_{2,-2,1}$ & $105 \left(15 s^2-14\right)
   \left(399375 s^{10}-1863400 s^8+3389440 s^6-3023632 s^4+1328128
   s^2-230400\right)$ \\
${}_{1,-1,1}$ & $-840 \left(5 s^2-4\right) \left(228125 s^{10}-549325
   s^8+255940 s^6+324664 s^4-352992 s^2+93824\right)$ \\
${}_{2,-1,1}$ & $210 \left(15
   s^2-14\right) \left(100875 s^{10}-275600 s^8+228220 s^6-6408 s^4-70272
   s^2+23296\right)$ \\
${}_{2,-1,2}$ & $-105 \left(5 s^2-4\right) \left(15 s^2-14\right)
   \left(13725 s^6-37680 s^4+34228 s^2-10304\right)$ \\
${}_{0,1,1}$ & $-1680 \left(5
   s^2-4\right)^2 \left(486525 s^8-1594290 s^6+1955772 s^4-1064576
   s^2+216960\right)$ \\
${}_{1,1,1}$ & $840 \left(5 s^2-4\right) \left(1531125 s^{10}-6503075
   s^8+10982780 s^6-9224760 s^4+3855648 s^2-641920\right)$ \\
${}_{2,1,1}$ & $-420 \left(15
   s^2-14\right) \left(61875 s^{10}-138825 s^8+51640 s^6+92200 s^4-89088
   s^2+22400\right)$ \\
${}_{1,1,2}$ & $1680 \left(5 s^2-4\right) \left(240750 s^8-775475
   s^6+932445 s^4-495822 s^2+98320\right)$ \\
${}_{2,1,2}$ & $105 \left(5 s^2-4\right)
   \left(226125 s^8-787950 s^6+1015020 s^4-572056 s^2+118944\right)$ \\
${}_{2,1,3}$ & $-840
   \left(15 s^2-14\right) \left(1125 s^6-3300 s^4+3235 s^2-1058\right)$ \\
${}_{0,2,1}$ & $-1680
   \left(5 s^2-4\right)^3 \left(41615 s^6-97838 s^4+76016 s^2-19488\right)$ \\
${}_{1,2,1}$ & $
   -1680 \left(5 s^2-4\right) \left(666875 s^{10}-2586600 s^8+4014940 s^6-3117320
   s^4+1210368 s^2-187904\right)$ \\
${}_{2,2,1}$ & $210 \left(5398125 s^{12}-27480750
   s^{10}+57999400 s^8-64973520 s^6+40757888 s^4-13579264 s^2+1878016\right)$ \\
${}_{1,2,2}$ & $
   -3360 \left(5 s^2-4\right)^2 \left(42850 s^6-108830 s^4+92099
   s^2-25958\right)$ \\
${}_{2,2,2}$ & $840 \left(5 s^2-4\right) \left(123750 s^8-359475
   s^6+378010 s^4-167724 s^2+25624\right)$ \\
${}_{2,2,3}$ & $-105 \left(639375 s^8-2259750
   s^6+2991200 s^4-1757840 s^2+387104\right)$ \\
${}_{0,3,1}$ & $-1120 \left(5 s^2-4\right)^2
   \left(270650 s^8-828285 s^6+945816 s^4-477232 s^2+89664\right)$ \\
${}_{1,3,1}$ & $280
   \left(5 s^2-4\right) \left(634500 s^{10}-2623725 s^8+4300340 s^6-3496152 s^4+1412352
   s^2-227584\right)$ \\
${}_{2,3,1}$ & $-70 \left(15 s^2-14\right) \left(76875 s^{10}-202950
   s^8+167700 s^6-16544 s^4-37376 s^2+12544\right)$ \\
${}_{0,3,2}$ & $-560 \left(5
   s^2-4\right)^2 \left(605775 s^6-1524950 s^4+1277728 s^2-356256\right)$ \\
${}_{1,3,2}$ & $560
   \left(5 s^2-4\right)^2 \left(9150 s^6-6435 s^4-9741 s^2+7154\right)$ \\
${}_{2,3,2}$ & $35
   \left(5 s^2-4\right) \left(104625 s^8-68850 s^6-286620 s^4+378936
   s^2-127232\right)$ \\
${}_{1,3,3}$ & $-280 \left(5 s^2-4\right) \left(52875 s^6-129225
   s^4+102900 s^2-26456\right)$ \\
${}_{2,3,3}$ & $-140 \left(15 s^2-14\right) \left(7875
   s^6-21225 s^4+19120 s^2-5756\right)$ \\
${}_{1,4,1}$ & $-840 \left(5 s^2-4\right)
   \left(516875 s^{10}-1956050 s^8+2950940 s^6-2217472 s^4+829440
   s^2-123392\right)$ \\
${}_{2,4,1}$ & $420 \left(5 s^2-4\right) \left(173625 s^{10}-668250
   s^8+1023720 s^6-780120 s^4+295872 s^2-44800\right)$ \\
${}_{0,4,2}$ & $-6720 \left(5
   s^2-4\right)^3 \left(730 s^4-1153 s^2+444\right)$ \\
${}_{1,4,2}$ & $-3360 \left(5
   s^2-4\right)^2 \left(66050 s^6-166215 s^4+139230 s^2-38808\right)$ \\
${}_{2,4,2}$ & $420
   \left(5 s^2-4\right)^2 \left(20925 s^6-38700 s^4+19984 s^2-1976\right)$ \\
${}_{1,4,3}$ & $
   840 \left(5 s^2-4\right) \left(10125 s^6-32150 s^4+33500 s^2-11456\right)$ \\
${}_{2,4,3}$ & $
   -2100 \left(5 s^2-4\right) \left(3525 s^6-9240 s^4+7996 s^2-2280\right)$ \\
${}_{1,5,1}$ & $
   -168 \left(5 s^2-4\right) \left(115000 s^{10}-434875 s^8+663700 s^6-512080 s^4+199936
   s^2-31488\right)$ \\
${}_{2,5,1}$ & $-105 s^2 \left(5 s^2-4\right) \left(15 s^2-14\right)
   \left(825 s^6-1990 s^4+1616 s^2-448\right)$ \\
${}_{0,5,2}$ & $-336 \left(5 s^2-4\right)^3
   \left(15425 s^4-24050 s^2+9112\right)$ \\
${}_{1,5,2}$ & $-84 \left(5 s^2-4\right)^2
   \left(551625 s^6-1390850 s^4+1167040 s^2-325728\right)$ \\
${}_{2,5,2}$ & $105 \left(5
   s^2-4\right)^2 \left(15 s^2-14\right) \left(225 s^4+288 s^2-364\right)$ \\
${}_{0,5,3}$ & $
   3360 \left(5 s^2-4\right)^2 \left(1575 s^4-2795 s^2+1256\right)$ \\
${}_{1,5,3}$ & $840
   \left(5 s^2-4\right) \left(8250 s^6-24975 s^4+25080 s^2-8336\right)$ \\
${}_{2,5,3}$ & $-420
   \left(5 s^2-4\right) \left(15 s^2-14\right)^2 \left(15 s^2-13\right)$ \\
${}_{2,6,1}$ & $35
   \left(5 s^2-4\right) \left(171375 s^{10}-616950 s^8+871680 s^6-600128 s^4+199168
   s^2-25088\right)$ \\
${}_{1,6,2}$ & $-280 \left(5 s^2-4\right)^3 \left(8265 s^4-12874
   s^2+4872\right)$ \\
${}_{2,6,2}$ & $-280 \left(5 s^2-4\right)^2 \left(11475 s^6-29280
   s^4+24828 s^2-6992\right)$ \\
${}_{0,6,3}$ & $560 \left(5 s^2-4\right)^3 \left(1335
   s^2-1166\right)$ \\
${}_{1,6,3}$ & $560 \left(5 s^2-4\right)^2 \left(8325 s^4-14910
   s^2+6764\right)$ \\
${}_{2,6,3}$ & $70 \left(5 s^2-4\right) \left(21375 s^6-61950 s^4+59960
   s^2-19328\right)$ \\
${}_{1,7,2}$ & $-60 \left(5 s^2-4\right)^3 \left(8385 s^4-13226
   s^2+5080\right)$ \\
${}_{0,7,3}$ & $10080 \left(5 s^2-4\right)^3 \left(90
   s^2-79\right)$ \\
${}_{1,7,3}$ & $12600 \left(5 s^2-4\right)^2 \left(105 s^4-191
   s^2+88\right)$ \\
${}_{2,8,2}$ & $-840 \left(3 s^2-2\right) \left(5 s^2-4\right)^3 \left(15
   s^2-14\right)$ \\
${}_{1,8,3}$ & $4200 \left(5 s^2-4\right)^3 \left(96
   s^2-85\right)$ \\
${}_{2,8,3}$ & $210 \left(5 s^2-4\right)^2 \left(525 s^4-990
   s^2+472\right)$ \\
${}_{1,9,3}$ & $7560 \left(5 s^2-4\right)^3 \left(10
   s^2-9\right)$ \\
${}_{2,10,3}$ & $315 \left(5 s^2-4\right)^3 \left(15 s^2-14\right)$ \\
${}_{0,0,3}$ & $2 \left(15 s^2-14\right)^3 \left(825 s^4-1445s^2+634\right)$ \\
${}_{0,0,2}$ & $-6 \left(5 s^2-4\right)^2 \left(2171250 s^8-7719525s^6+10225470 s^4-5983260 s^2+1305248\right)$ \\
${}_{1,0,2}$ & $-3 \left(5 s^2-4\right)
   \left(15 s^2-14\right)^2 \left(1800 s^6+2655 s^4-8208 s^2+3928\right)$ \\
${}_{0,0,1}$ & $48
   \left(5 s^2-4\right)^2 \left(9060750 s^{10}-34431275 s^8+51858720 s^6-38675200
   s^4+14258176 s^2-2072064\right)$ \\
${}_{1,0,1}$ & $-12 \left(5 s^2-4\right) \left(93223125
   s^{12}-421210500 s^{10}+784654200 s^8-771469840 s^6+422629664 s^4-122600960
   s^2+14780416\right)$ \\
${}_{2,0,1}$ & $6 \left(15 s^2-14\right) \left(2328750 s^{12}-8703375
   s^{10}+13317150 s^8-10848180 s^6+5157560 s^4-1450624 s^2+200704\right)$ \\
\noalign{\smallskip}\hline
\end{tabular}\normalsize
\end{table*}

The computation of additional orders finds similar structures. Thus, for instance, the fourth order term of the new Hamiltonian takes the form
\begin{equation} \label{K04g}
K_{0,4}=\epsilon^4\frac{\mu}{r}\frac{p^2}{r^2}\frac{9s^2}{1280(5s^2-4)^6}\sum_{j=0}^5\frac{p^j}{r^j}\sum_{k=0}^{\lfloor3-\frac{1}{2}j\rfloor}e^{2k}\gamma_{4,j,k},
\end{equation}
and the fourth order term of the generating function takes the form
\begin{eqnarray} \nonumber
W_4 &=& \frac{\epsilon^4G}{30105600(5s^2-4)^6}\sum_{l=1}^4\sum_{\substack{k=l-6 \\ k\ne0}}^{2l+6}\sum_{j=0}^6\Gamma_{4,j,k,l}e^{2j+k^*} \\ \label{W4g}
&& \times s^{2l}\sin(kf+2l\omega)+C_4.
\end{eqnarray}
While the former involves 15 inclination polynomials $\gamma_{4,j,k}$, the later comprises up to 124 non-vanishing inclination polynomials $\Gamma_{4,j,k,l}$, of degree 9 in $s^2$. Therefore, these and other polynomials resulting from following orders are not listed due to their length.
\par

The normalization of $G$ has been extended up to the order six without major trouble, except for the increase of the computational burden of successive orders due to the notable growth of the length of the series to be handled, on the one hand, and the increasing size of the rational coefficients resulting from the integer arithmetic used, on the other. Thus, the new Hamiltonian terms $K_{0,5}$ and $K_{0,6}$ take analogous forms to the previous orders, with 24 inclination polynomials $\gamma_{5,j,k}$, and 35 $\gamma_{6,j,k}$. Similarly, the generating function terms $W_{5}$ and $W_{6}$ have been arranged in the same form as previous orders of the perturbation approach, with 254 non-vanishing coefficients $\Gamma_{5,j,k,l}$ and 429 $\Gamma_{6,j,k,l}$, 15 of which correspond to the integration constant $C_5$ and 21 to $C_6$.

At the end, since the generating function is known, the transformation equations from the original Delaunay variables to the new ones, and vice versa, are readily computed by standard application of the fundamental recursion (\ref{Lie:triangle}) (see \cite{Deprit1969} for details). The procedure ends by replacing the original variables by the new ones in the computed terms $K_{0,m}$.
\par

We remark that the procedure described here is not equivalent to the combination of the elimination of the parallax  and the elimination of the perigee into a single transformation, in spite of the fact that the total angular momentum is converted into a formal integral in both cases. On the contrary, the perigee has been removed here keeping as much short-period terms as possible in the new Hamiltonian. With this strategy, the size of the generating function of the first partial normalization is astonishingly smaller than the one that would be obtained with other alternatives in the literature. To check that, we fully expanded both the Hamiltonian and the generating function and reckoned the number of terms. This procedure is the one that has been traditionally used in the literature to assess the complexity of a perturbation solution \cite{DepritRom1970,CoffeyDeprit1982}. Thus, for instance, we find 5, 56, 367, 1152, 2627, and 4897 expanded terms for the 1st, 2nd \dots, 6th order terms of the generating function, respectively, whereas much longer expressions have been reported in the literature. For instance, the sixth order of the generating function reported in  \cite{SanJuanetal2013} entails 39630 terms, which is almost one order of magnitude larger than the one obtained with the current approach. It is worth mentioning that further simplifications could be obtained in particular cases, like when constraining the application of the analytical solution to the case of low eccentricity orbits, in which case many of the terms involved can be neglected \cite{Lara2008,GaiasColomboLara2020}.
\par

The comparisons are, nevertheless, inconclusive due to the diversity of approaches used in the computation of the variety of solutions reported in the literature, which involve representation in different variables, on the one hand, and yield distinct one-degree-of-freedom Hamiltonians, on the other. The partially normalized Hamiltonian obtained here is definitely longer than the one obtained after the classical elimination of the perigee. However, this is not of concern in the complete, as opposite to partial, normalization of the main problem. Indeed, after the following elimination of short-period terms either procedure should arrive to the same Hamiltonian. We didn't find reported data for the second normalization after the elimination of the perigee, but one would expect that the figures should be balanced to some extent ---the generating function of the short-period elimination probably being heavier in our case than in other prospective approaches. While the needed data for the thorough comparison is not available, our claims must restrict to the checked fact that our approach makes both transformations of manageable size ---as we will show in the next section, where the short-period effects are removed in a Delaunay normalization. This feature makes the evaluation of higher orders of our perturbation solution certainly practicable.

\subsection{Delaunay normalization}

The partially normalized Hamiltonian is just of one degree of freedom, yet it is not separable. To get an explicit analytical solution we remove the short-period terms that are associated to the radius by means of a Delaunay normalization \cite{Deprit1982}. The partially reduced Hamiltonian from which we start takes again the form of Eq.~(\ref{mainP}), but it is now written in prime Delaunay variables $(\ell',g',h',L',G',H')$. In this Hamiltonian the variables $g'$ and $h'$ are cyclic, and hence $H'=H$ and $G'$ remain constant for given initial conditions. The zeroth order term is the same as in Eq.~(\ref{mainPK0}), the term $K_{1,0}$ is given by Eq.~(\ref{K01g}), whereas terms $K_{m,0}$ with $m\ge2$ are no longer void, and, on the contrary are given by Eqs.~(\ref{K02g}), (\ref{K03g}), (\ref{K04g}), for $m=2,3,4$, respectively, and analogous equations for higher orders that, due to its length, are not printed in the paper. Moreover, the definition in Eq.~(\ref{mieps}) turns $\epsilon$ into a physical constant, rather than a function, when witten in the prime variables. In consequence $\epsilon$ can replace now the formal small parameter $\varepsilon$ used before in the Lie transforms procedure, assumed, of course, that the corresponding scaling of the Hamiltonian terms is properly made.
\par

The homological equation to be solved at each step of the new Lie transformation is the same as before, either in the form given by Eq.~(\ref{homosoll}) or Eq.~(\ref{homosol}), except for choosing a non-zero integration constant does not provide any advantage in the current case. In fact, since the new Hamiltonian must be free from terms depending on the mean anomaly, which are obtained by averaging, the homological equation can be further particularized. That is,
\begin{equation} \label{averaging}
K_{0,m}=\langle\widetilde{K}_{0,m}\rangle_\ell=\frac{1}{2\pi}\int_0^{2\pi}\widetilde{K}_{0,m}\frac{r^2}{a^2\eta}\,\mathrm{d}f,
\end{equation}
and hence
\begin{eqnarray} \nonumber
W_m &=& -\frac{\ell}{n}K_{0,m}+\frac{1}{n}\int\widetilde{K}_{0,m}\frac{r^2}{a^2\eta}\,\mathrm{d}f \\
\label{homolo}
&=& K_{0,m}\frac{\phi}{n}+\frac{1}{n}\int\left(\widetilde{K}_{0,m}\frac{r^2}{a^2\eta}-K_{0,m}\right)\mathrm{d}f,
\end{eqnarray}
in which $\phi=f-\ell$ is the equation of the center, and the terms under the integral sign in the final form of the equation are purely periodic in $f$.
\par

Thus, at the first order of the Lie transforms procedure we find
\begin{equation} \label{K01l}
K_{0,1}=\langle{K}_{1,0}\rangle_\ell\equiv \frac{\mu}{p}\eta^3 \left(3 s^2-2\right),
\end{equation}
from which
\begin{equation} \label{W1l}
W_1=G'\left(3s^2-2\right)(e\sin{f}+\phi ).
\end{equation}
Recall that the symbols $p$, $\eta$, $e$, etc. are now functions of the Delaunay prime variables.
\par

At the second order, the known terms are given again by Eqs.~(\ref{triangle2}) and (\ref{K11}), from which, using Eq.~(\ref{averaging}),
\begin{equation} \label{K02l}
K_{0,2}=-\frac{\mu}{p}\frac{3}{4}\eta^3\sum_{j=0}^2\lambda_{2,j}\eta^j,
\end{equation}
with $\lambda_{2,0}=5(7s^4-16s^2+8)$, $\lambda_{2,1}=4(3s^2-2)^2$, and $\lambda_{2,2}=5s^4+8s^2-8$. The second order term of the generating function is then trivially integrated from Eq.~(\ref{homolo}), to yield
\begin{eqnarray} \nonumber
W_2 &=& 
-\frac{G'\beta}{32(5s^2-4)^2}\sum_{j=1}^3\sum_{k=0}^{3-\lfloor\frac{1}{2}j\rfloor}\Lambda_{2,j,k}\eta^ke^j\sin{jf} \\ \label{W2l}
&& -G'\frac{3}{4}\phi\sum_{j=0}^1\Phi_{2,j}e^{2j},
\end{eqnarray}
where $\beta=1/(1+\eta)$, $\Phi_{2,0}=8(s^2-1)(5s^2-4)$, $\Phi_{2,1}=8-8s^2-5s^4$, and the inclination polynomials $\Lambda_{2,j,k}$ are provided in Table~\ref{t:Lambda2}.

\begin{table}[htb]
\caption{Inclination polynomials $\Lambda_{2,j,k}$ in Eq.~(\protect\ref{W2l}). }
\label{t:Lambda2}
\begin{tabular}{@{}ll@{}}
\hline\noalign{\smallskip}
${}_{1,0}$ & $15 (3 s^2-2) (805 s^6-2448 s^4+2400 s^2-768)$ \\
${}_{1,1}$ & $3 (3s^2-2) (2225 s^6-8160 s^4+8928 s^2-3072)$ \\
${}_{1,2}$ & $3 (825 s^8-3030s^6+4064 s^4-2368 s^2+512)$ \\
${}_{1,3}$ & $-3 s^2 (975 s^6-2250 s^4+1728 s^2-448)$ \\
${}_{2,0}$ & $6(1925 s^8-6210 s^6+7452 s^4-3936 s^2+768)$ \\
${}_{2,1}$ & $6 (125 s^8-930 s^6+1660s^4-1120 s^2+256)$ \\
${}_{3,0}$ & $2625 s^8-7270 s^6+7408 s^4-3264 s^2+512$ \\
${}_{3,1}$ & $s^2 (825 s^6-1990s^4+1616 s^2-448)$ \\
\noalign{\smallskip}\hline
\end{tabular}
\end{table}

Integrals cease to be trivial at the third order. Indeed, the formal computation of $K_{0,3}$ using the fundamental recursion (\ref{Lie:triangle}) yields two different types of terms.
\par

The first type consists of terms depending on the equation of the center, which can be reduced to the form
\[
\frac{\mu}{r}\frac{p}{r}P(e)Q(s)\phi\sin{mf},
\]
with $P$ and $Q$ denoting arbitrary eccentricity and inclination polynomials. Definite integration of these kinds of terms is carried out from expressions in \cite{Metris1991}, whereas the indefinite integration is achieved by parts, to get
\[
ma^2\eta\int \frac{\sin{mf}}{r^2}\phi\,\mathrm{d}\ell=\frac{\sin{mf}}{m}-\phi\cos{mf}-\int\cos{mf}\mathrm{d}\ell,
\]
in which the antiderivatives of cosines of multiples of the true anomaly are carried out after expressing them in terms of $r$ and $R=\mathrm{d}r/\mathrm{d}t$, rather than in trigonometric functions, as discussed in \cite{Jefferys1971}.
\par

The second type consists of terms free from the equation of the center. In these terms, the trigonometric functions of the true anomaly can be replaced by inverse powers of the radius, without involving the radial velocity. We found that the exponents of the inverse of $r$ range from 0 to 8 missing the exponent 1. Therefore, both definite and indefinite integrals of terms of the second type are trivially solved.
\par

Proceeding in this way, we compute the third order Hamiltonian term
\begin{equation} \label{K03l}
K_{0,3}=\frac{\mu}{p}\frac{9\eta^3}{16(5s^2-4)^2}\sum_{j=0}^4\lambda_{3,j}\eta^j,
\end{equation}
with the inclination polynomials $\lambda_{3,j}$ given in Table~\ref{t:lambda3}.
\par

\begin{table}[htb]
\caption{Inclination polynomials $\lambda_{3,j}$ in Eq.~(\protect\ref{K03l}).}
\label{t:lambda3}
\begin{tabular}{@{}ll@{}}
\hline\noalign{\smallskip}
${}_{0}$ & $5(28700 s^{10}-107205 s^8+158960 s^6-118492 s^4$ \\
 & $+45152 s^2-7168)$ \\
${}_{1}$ & $60(3 s^2-2)(5 s^2-4)^2(7 s^4-16 s^2+8)$ \\
${}_{2}$ & $-2(28675s^{10}-98005 s^8+130852 s^6-87164 s^4$ \\
 & $+30176 s^2-4608)$ \\
${}_{3}$ & $20(3 s^2-2)(5 s^2-4)^2(5 s^4+8 s^2-8)$ \\
${}_{4}$ & $-s^2(15 s^2-14)(450s^6-925 s^4+590 s^2-112)$ \\
\noalign{\smallskip}\hline
\end{tabular}
\end{table}

The corresponding term of the generating function is
\begin{eqnarray} \nonumber
W_3 &=& \frac{G'\beta^2}{128(5s^2-4)^3}\sum_{j=1}^6\sum_{k=0}^{7-2\lfloor\frac{1}{2}j\rfloor}\Lambda_{3,j,k}\eta^{k-1}e^j\sin{jf} \\ \label{W03l}
&& +\frac{3G'}{16(5s^2-4)^2}\phi\sum_{j=0}^3\sum_{k=0}^{4}\Phi_{3,j,k}\eta^ke^j\cos{jf},
\end{eqnarray}
in which the inclination polynomials $\Phi_{3,j,k}$ are given in Table \ref{t:W3lc}, and those $\Lambda_{3,j,k}$ in Table \ref{t:W3ls}.

\begin{table}[htb]
\caption{Non-zero coefficients $\Phi_{3,j,k}$ in Eq.~(\protect\ref{W03l}). $\Phi_{3,1,0}=\frac{15}{2}\Phi_{3,0,3}$, $\Phi_{3,1,2}=-\frac{3}{2}\Phi_{3,0,3}$, $\Phi_{3,2,0}=3\Phi_{3,0,3}$, and $\Phi_{3,3,0}=\frac{1}{2}\Phi_{3,0,3}$. }
\label{t:W3lc}
\begin{tabular}{@{}ll@{}}
\hline\noalign{\smallskip}
${}_{0,0}$ & $5 (89100 s^{10}-323615 s^8+466320 s^6-337684 s^4$ \\
 & $+125216 s^2-19456)$ \\
${}_{0,2}$ & $-2(112125 s^{10}-374775 s^8+488460 s^6-314932 s^4$ \\
 & $+103840 s^2-14848)$ \\
${}_{0,3}$ & $8 (3 s^2-2) (5 s^2-4)^2 (5 s^4+8 s^2-8)$ \\
${}_{0,4}$ & $-3 s^2 (15s^2-14) (450 s^6-925 s^4+590 s^2-112)$ \\
\noalign{\smallskip}\hline
\end{tabular}
\end{table}
\begin{table*}[htbp]
\caption{Non-zero inclination polynomials $\Lambda_{3,j,k}$ in Eq.~(\protect\ref{W03l}). }
\label{t:W3ls}
\begin{tabular}{@{}ll@{}}
\hline\noalign{\smallskip}
${}_{1,0}$ & $-864 (3 s^2-2)^3 (5 s^2-4)^3$ \\
${}_{1,1}$ & $3 (22218875 s^{12}-104346550 s^{10}+202703740 s^8-209869352s^6+123038240s^4-39033472 s^2+5275648)$ \\
${}_{1,2}$ & $12(10925500s^{12}-50711075s^{10}+97386820s^8-99715748s^6+57863024s^4-18199872s^2+2445312)$ \\
${}_{1,3}$ & $3(27560125s^{12}-119080550s^{10}+212650740s^8-202245448s^6+109190304s^4-32208768s^2+4128768)$ \\
${}_{1,4}$ & $12(3155125s^{12}-10820800s^{10}+13899620s^8-7620256s^6+944256s^4+613248s^2-167936)$ \\
${}_{1,5}$ & $3(5410625s^{12}-17331450s^{10}+19448180s^8-6842968s^6-2742560s^4+2556544s^2-491520)$ \\
${}_{1,6}$ & $-12(59625s^{12}-415275s^{10}+942920s^8-994980s^6+529776s^4-134592s^2+12288)$ \\
${}_{1,7}$ & $-3 s^2 (77625 s^{10}-568950 s^8+1256420s^6-1222216 s^4+550816 s^2-94080)$ \\ 
${}_{2,0}$ & $-1044 (3 s^2-2)^3(5 s^2-4)^3$ \\
${}_{2,1}$ & $24 (131000s^{12}-1121875s^{10}+3061340s^8-3989664s^6+2758768s^4-983648s^2+143360)$ \\
${}_{2,2}$ & $96(51625s^{12}-437800s^{10}+1183290s^8-1528682s^6+1049344s^4-372240s^2+54144)$ \\
${}_{2,3}$ & $-12(5s^2-4)(16375s^{10}+64070s^8-257508 s^6+297320 s^4-145792 s^2+26624)$ \\
${}_{2,4}$ & $-12(263625s^{12}-1000750s^{10}+1526820s^8-1206712s^6+539616s^4-142592s^2+19968)$ \\
${}_{2,5}$ & $-12s^2 (162375 s^{10}-576100 s^8+787020s^6-506248 s^4+145984 s^2-12992)$ \\
${}_{3,0}$ & $-656 (3 s^2-2)^3(5 s^2-4)^3$ \\
${}_{3,1}$ & $-934875s^{12}+605000s^{10}+4973120s^8-10412952s^6+8554272s^4-3281280s^2+491520$ \\
${}_{3,2}$ & $-2080125s^{12}+3562000s^{10}+2881300s^8-11103360s^6+10155456s^4-4038912s^2+614400$ \\
${}_{3,3}$ & $-(5 s^2-4) (254925 s^{10}-526480 s^8+318084 s^6-10672s^4-40064 s^2+8192)$ \\
${}_{3,4}$ & $3(28875s^{12}-147800s^{10}+254260s^8-160992s^6-11968s^4+54016s^2-16384)$ \\
${}_{3,5}$ & $-12 s^2 (4500s^{10}-4125 s^8-16075 s^6+31670 s^4-20664 s^2+4704)$ \\
${}_{4,0}$ & $-240 (3s^2-2)^3 (5 s^2-4)^3$ \\
${}_{4,1}$ & $6 (5 s^2-4)(50325 s^{10}-157660 s^8+180520 s^6-90312 s^4+18112 s^2-1024)$ \\
${}_{4,2}$ & $12 (5 s^2-4) (47475 s^{10}-147000 s^8+164852 s^6-79056s^4+14208 s^2-512)$ \\
${}_{4,3}$ & $6 s^2 (5 s^2-4) (44625s^8-136340 s^6+149184 s^4-67800 s^2+10304)$ \\
${}_{5,0}$ & $-48 (3s^2-2)^3 (5 s^2-4)^3$ \\
${}_{5,1}$ & $6 s^2 (5 s^2-4)^2(180 s^6-609 s^4+530 s^2-112)$ \\
${}_{5,2}$ & $3 s^2 (5 s^2-4)(1125 s^8-7440 s^6+11516 s^4-6144 s^2+896)$ \\
${}_{5,3}$ & $-3 s^4 (5s^2-4) (15 s^2-14) (45 s^4+36 s^2-56)$ \\
${}_{6,0}$ & $-4(3 s^2-2)^3 (5 s^2-4)^3$\\
\noalign{\smallskip}\hline
\end{tabular}
\end{table*}

In the process of carrying out the normalization to higher orders we need to deal with trigonometric series of notably increasing length. However, we only found terms of the same two types as before in the solution of the homological equation, which, in consequence, is analogously integrated. In this way, the second normalization has been extended up to the order 6 in the small parameter, in agreement with the order to which the perigee was previously eliminated. Corresponding Hamiltonian and generating function terms are analogously arranged in the form of Eqs.~(\ref{K03l}) and (\ref{W03l}), respectively, yet the inclination polynomials comprise much more longer listings, as expected. Thus, for instance, at the fourth order the normalized Hamiltonian term
\begin{equation} \label{K04l}
K_{0,4}=\frac{\mu}{p}\frac{9\eta^3}{64(5s^2-4)^3}\sum_{j=0}^6\lambda_{4,j}\eta^j,
\end{equation}
contributes 7 new inclination polynomials, which are polynomials of degree 7 in the square of the sine of the inclination, whereas the generating function term
\begin{eqnarray} \nonumber
W_4 &=& 
\frac{G'\beta^3}{20480(5s^2-4)^6}\sum_{j=1}^6\sum_{k=0}^{7-2\lfloor\frac{1}{2}j\rfloor}\Lambda_{4,j,k}\eta^{k-3}e^j\sin{jf} \\ \label{W04l}
&& +\frac{3G'}{256(5s^2-4)^3}\phi\sum_{j=0}^5\sum_{k=0}^{6}\Phi_{4,j,k}\eta^ke^j\cos{jf},
\end{eqnarray}
is made of 20 non-zero coefficients $\Phi_{4,j,k}$, of degree 7 in $s^2$, and 72 non-zero coefficients $\Lambda_{4,j,k}$, which are of degree 10 in $s^2$. Note that, because of the denominators factoring the summations, the maximum degree amounts to 4 in any case, in agreement with the order of the perturbation. At the 5th order we find 9 coefficients of the type $\lambda$ and 41 nonvanishing coefficients of the type $\Phi$, which are of degree 11 in $s^2$ although they are divided by $(5s^2-4)^6$, as well as 134 nonvanishing trigonometric polynomials of the type $\Lambda$, which are of degree 12 in $s^2$ yet they must be divided by $(5s^2-4)^7$. Finally, the figures of the 6th order are 11 and 70 of degree 13 in $s^2$ for $\lambda$ and $\Phi$, respectively, and 218 for $\Lambda$ of degree 16.
\par

Once reached the desired order of the second normalization, the procedure ends by changing prime variables by double-prime variables in the new Hamiltonian terms $K_{0,m}$.
\par

In order to provide comparative figures of the computational burden of this normalization with other approaches that might be carried out, we expanded the series that comprise the solution and reckon the number of separate terms, as we already did in the normalization of the total angular momentum. We found that the 1st, 2nd, \dots 6th-order terms of the generating function of the Delaunay normalization comprise 4, 48, 257, 931, 2266, and 4826 terms respectively. Note that different arrangements from the one chosen by us ---Eqs.~(\ref{W03l}), (\ref{W04l}), etc.--- may provide different figures than those reported, yet should be of analogous magnitude. Following the same procedure with the completely reduced Hamiltonian, we reckon up to 2, 9, 29, 55, 106, and 152 coefficients, for the 1st, 2nd, \dots 6th-order term, respectively.

\subsection{Secular terms}

After neglecting higher order effects of $J_2$, the normalized Hamiltonian is
\begin{equation} \label{Hn}
\mathcal{K}''=\mathcal{K}''(-,-,-,L'',G'',H'')\equiv\sum_{m\ge0}^{\tilde{m}}\frac{\epsilon^m}{m!}K_{m,0},
\end{equation}
in which $\tilde{m}\le6$ for the different approximations provided by the computed perturbation solution, and terms $K_{m,0}$ are obtained by replacing prime by double-prime variables in corresponding terms given by Eqs.~(\ref{K01l}), (\ref{K02l}), (\ref{K03l}), as well as higher order Hamiltonian terms $K_{0,m}$ that have not been displayed. That is, the symbols $p$, $\eta$, and $s$ in these equations are assumed to be functions of the double-prime Delaunay momenta. Recall that $\epsilon$ is obtained from Eq.~(\ref{mieps}) by making $p=G''^2/\mu$.

The corresponding solution to the flow stemming from Eq.~(\ref{Hn}) is
\begin{eqnarray*}
\ell'' &=& \ell''_0+n_\ell{t}, \\
g''&=& g''_0+{n}_\omega{t} \\
h'' &=& h''_0+n_ht,
\end{eqnarray*}
in which, from Hamilton equations,
\[
n_\ell=\frac{\partial\mathcal{K}''}{\partial{L}''}, \qquad
n_g=\frac{\partial\mathcal{K}''}{\partial{G}''}, \qquad
n_h=\frac{\partial\mathcal{K}''}{\partial{H}''},
\]
$L''=L''_0$, $G''=G'_0$, and $H''=H_0$, are the initialization constants of the analytical solution, and $\ell''_0$, $g''_0$, $h''_0$, $L''_0$, and $G'_0$, are computed applying consecutively the inverse transformation of the normalization of the angular momentum and the Delaunay normalization to corresponding initial conditions of the original problem.
\par

On the other hand, Delaunay variables are singular for equatorial orbits, in which the argument of the node is not defined, as well as for circular orbits, in which the argument of the periapsis is not defined. In particular, the singularity for circular orbits reveals immediately by the appearance of the eccentricity in denominators of the transformation equations of both the mean anomaly and the argument of the perigee (see Eqs.~(20) and (21) of \cite{Brouwer1959}, for instance). This fact not only makes singular the transformation of these elements for exactly circular orbits, but prevents convergence of the respective perturbation series for small values of the eccentricity.
\par

Different sets of non-singular variables can be used to avoid these issues \cite{Lyddane1963}. In particular, troubles related with low eccentricities are commonly avoided by replacing the mean anomaly, the argument of the periapsis, and the total angular momentum with the non-canonical variables given by the mean distance to the node $F=\ell+g$, also called mean argument of the latitude, and the components of the eccentricity vector in the nodal frame $C=e\cos\omega$, $S=e\sin\omega$, also called semi-equinoctial elements \cite{Konopliv1990,Cook1992}. In these variables, the secular terms of the main problem are given by \cite{DepritRom1970}
\begin{eqnarray} \label{Fsecular}
F &=& F_0+n_Ft, \\ \label{Csecular}
C &=& e\cos(g_0+n_g{t})=C_0\cos{n}_g{t}-S_0\sin{n}_g{t}, \\ \label{Ssecular}
S &=& e\sin(g_0+n_g{t})=S_0\cos{n}_g{t}+C_0\sin{n}_g{t}, \\ \label{Lsecular}
L &=& L_0, \\ \label{hsecular}
h &=& h_0+n_ht, \\ \label{Hsecular}
H &=& H_0,
\end{eqnarray}
in which $n_F=n_\ell+n_g$, $C_0=e\cos{g}_0$, $S_0=e\sin{g}_0$, and $e=(1-G_0^2/L_0^2)^{1/2}$. The double-prime notation has been omitted for simplicity.
\par

After standard partial differentiation, we obtain
\begin{eqnarray} \label{nFsecular}
n_F &=& n+n\sum_{m=1}^{\tilde{m}}
\frac{\epsilon^m}{(5s^2-4)^m}\,\sum_{i=0}^{2m-1}\Psi_{m,i}(s)\eta^i, \\ \label{ngsecular}
n_g &=& n\sum_{m=1}^{\tilde{m}}
\frac{\epsilon^m}{(5s^2-4)^m}\sum_{i=0}^{2m-2}\omega_{m,i}(s)\eta^i, \\ \label{nhsecular}
n_h &=& n{c}\sum_{m=1}^{\tilde{m}}
\frac{\epsilon^m}{(5s^2-4)^m}\sum_{i=0}^{2m-2}\Omega_{m,i}(s)\eta^i,
\end{eqnarray}
where the inclination polynomials $\Psi_{m,i}$, $\omega_{m,i}$, and $\Omega_{m,i}$ are provided in Tables \ref{t:Psi}, \ref{t:omega}, and \ref{t:Omega}, respectively, up to the third order of the perturbation approach. In these tables we can check that, in spite of the general arrangement of the secular frequencies in Eqs.~(\ref{nFsecular})--(\ref{nhsecular}), the critical inclination divisors $5s^2-4$ start to appear only at the third order truncation.
\par

\begin{table}[htbp]
\caption{Inclination polynomials $\Psi_{i,j}$ in Eq.~(\protect\ref{nFsecular}).}
\label{t:Psi}
\begin{tabular}{@{}llll@{}}
\hline\noalign{\smallskip}
${}_{1,0}$ & $-3 (5 s^2-4)^2$ \\
${}_{1,1}$ & $-3 (3 s^2-2) (5s^2-4)$ \\
${}_{2,0}$ & $\frac{15}{8} (5 s^2-4)^2 (77 s^4-172s^2+88)$ \\
${}_{2,1}$ & $\frac{9}{8} (5 s^2-4)^2 (155 s^4-256s^2+104)$ \\
${}_{2,2}$ & $\frac{3}{8} (5 s^2-4)^2 (189 s^4-156s^2+8)$ \\
${}_{2,3}$ & $\frac{15}{8} (5 s^2-4)^2 (5 s^4+8s^2-8)$ \\
${}_{3,0}$ &$-\frac{15}{32}(2439500s^{12}-11312175s^{10}+21772080s^8$ \\
 & $-22346500s^6+12956400s^4-4043136 s^2+533248)$ \\
${}_{3,1}$ & $-\frac{45}{32}(5 s^2-4) (62300 s^{10}-260365 s^8+431504 s^6$ \\
 & $-356508s^4+147552s^2-24576)$ \\
${}_{3,2}$ & $\frac{3}{16} (1835625 s^{12}-7723875 s^{10}+13291500s^8$ \\
 & $-12015300 s^6+6064176 s^4-1644928 s^2+192256)$ \\
${}_{3,3}$ & $\frac{15}{16}(5s^2-4)(18175s^{10}-85105 s^8+153172 s^6-136540 s^4$ \\
 & $+61408s^2-11264)$ \\
${}_{3,4}$ & $\frac{3}{32} (213750 s^{12}-1441125 s^{10}+3537000s^8-4313100 s^6$ \\
 & $+2835280 s^4-967808 s^2+135424)$ \\
${}_{3,5}$ & \multicolumn{3}{l}{$\frac{21}{32} s^2(5 s^2-4) (15 s^2-14)(450 s^6-925 s^4+590s^2-112)$}\\
\noalign{\smallskip}\hline
\end{tabular}
\end{table}
\begin{table}[htbp]
\caption{Inclination polynomials $\omega_{i,j}$ in Eq.~(\protect\ref{ngsecular}). }
\label{t:omega}
\begin{tabular}{@{}ll@{}}
\hline\noalign{\smallskip}
${}_{1,0}$ & $-3 (5 s^2-4)^2$ \\
${}_{2,0}$ & $\frac{15}{8} (5s^2-4)^2 (77 s^4-172 s^2+88)$ \\
${}_{2,1}$ & $9 (3 s^2-2)(5 s^2-4)^3$ \\
${}_{2,2}$ & $\frac{3}{8} (5 s^2-4)^2 (45 s^4+36s^2-56)$ \\
${}_{3,0}$ & $-\frac{15}{32}(2439500s^{12}-11312175s^{10}+21772080s^8$ \\
 & $-22346500s^6+12956400s^4-4043136 s^2+533248)$ \\
${}_{3,1}$ & $-\frac{45}{4} (5 s^2-4)^3 (168 s^6-497 s^4+460s^2-136)$ \\
${}_{3,2}$ & $\frac{3}{16} (2150625 s^{12}-9409875 s^{10}+16968300s^8$ \\
 & $-16218180 s^6+8729136 s^4-2535808 s^2+315136)$ \\
${}_{3,3}$ & $-\frac{15}{4}(5 s^2-4)^3 (105 s^6+39 s^4-228 s^2+104)$ \\
${}_{3,4}$ & $\frac{3}{32} (438750 s^{12}-1771125 s^{10}+2865000 s^8-2345100 s^6$ \\
 & $+999760s^4-199808 s^2+12544)$ \\
\noalign{\smallskip}\hline
\end{tabular}
\end{table}
\begin{table}[htb]
\caption{Inclination polynomials $\Omega_{i,j}$ in Eq.~(\protect\ref{nhsecular}). }
\label{t:Omega}
\begin{tabular}{@{}ll@{}}
\hline\noalign{\smallskip}
${}_{1,0}$ & $-6 (5 s^2-4)$ \\
${}_{2,0}$ & $\frac{15}{2} (5s^2-4)^2 (7 s^2-8)$ \\
${}_{2,1}$ & $18 (3 s^2-2) (5s^2-4)^2$ \\
${}_{2,2}$ & $\frac{3}{2} (5 s^2-4)^2 (5s^2+4)$ \\
${}_{3,0}$ & $-\frac{15}{8}(215250s^{10}-823025 s^8+1255040s^6-953760s^4$ \\
 & $+361088s^2-54464)$ \\
${}_{3,1}$ & $-\frac{45}{4}(5s^2-4)^3 (63 s^4-124 s^2+56)$ \\
${}_{3,2}$ & $\frac{3}{8}(430125s^{10}-1553550s^8+2222340s^6-1570224s^4$ \\
 & $+546432s^2-74624)$ \\
${}_{3,3}$ & $-\frac{15}{4} (5 s^2-4)^3 (45 s^4+28 s^2-40)$ \\
${}_{3,4}$ & $\frac{3}{8}(50625s^{10}-168375s^8+215900s^6-130800s^4$ \\
 & $+35840s^2-3136)$ \\
\noalign{\smallskip}\hline
\end{tabular}
\end{table}

Computing the periodic corrections of the semi-equi\-noctial variables would require to carry out expansions of the true anomaly. Therefore, we rather formulate them in polar variables. In this way we avoid the trouble of small denominators in the $J_2$-problem, yet, of course, the nodes of exactly equatorial orbits remain undefined. For simplicity we do not deal with this latter case, which, if desired, could be approached using different sets of non-singular variables.
\par

\section{Performance of the solution}

The performance of the analytical solution is illustrated in three different cases. The first one is a low-altitude, almost-circular orbit with the orbital parameters of the PRISMA mission \cite{PerssonJacobssonGill2005}. The second is a high-eccentricity orbit with the typical characteristics of a geostationary transfer orbit (GTO), which we borrowed from \cite{Konopliv1991}. The third test has been specifically carried out to check the behavior of the analytical solution when approaching the critical inclination resonance. For this last case we selected an orbit with orbital parameters similar to the TOPEX orbit, which departs only $\sim3^\circ$ from the inclination resonance condition \cite{Frauenholzetal1998}. Orbital elements corresponding to these three cases are presented in Table \ref{t:testcases}. The reference, true orbits have been propagated numerically in extended precision to assure that all the computed points are accurate with 15 significant digits in the decimal representation. We requested that precision because this is the maximum number of digits with which exact decimal operations are guaranteed in standard double-precision \cite{Kahan1997}. Because the analytical solution is evaluated in double-precision the reference numeric orbit is considered exact.
\par

\begin{table}[htb] \tabcolsep 5pt
\caption{Initial conditions of the test cases. Angles are in degrees.}
\label{t:testcases}
\begin{tabular}{@{}lllllrr@{}} \multicolumn{1}{c}{type}
 & \multicolumn{1}{c}{$a$ (km)} & \multicolumn{1}{c}{$e$} & \multicolumn{1}{c}{$I$}
 & \multicolumn{1}{c}{$\Omega$} & \multicolumn{1}{c}{$\omega$} & \multicolumn{1}{c}{$\ell$} \\[0.33ex]
\hline\noalign{\smallskip}
PRISMA & $6878.137$ & $0.001$  & $97.42$ & $168.162$ &  $20$ &  $30$ \\
TOPEX  & $7707.270$ & $0.0001$ & $66.04$ & $180.001$ & $270$ & $180$ \\
GTO    & $24460.00$ & $0.73$   & $30.$   & $170.1$   & $280$ &   $0$ \\
\noalign{\smallskip}\hline
\end{tabular}
\end{table}

Two different kind of errors are associated to the truncation order of the analytical solution. On the one hand, the truncation of the secular terms introduces an error in the computation of the frequencies of the analytical solution, Eqs.~(\ref{nFsecular})--(\ref{nhsecular}), which will make the perturbation solution to degrade with time. On the other hand, the truncation affects also the generating function from which the  periodic corrections are derived. The latter fundamentally affects the amplitude of the periodic errors, and, therefore, the quality of the ephemeris provided by the analytical solution is not expected to deteriorate significantly with time by effect of this error. However, this is only true for the \emph{direct} transformation, a case in which we will see that it is acceptable to use a lower order truncation than the truncation used for the secular terms. On the contrary, the truncation order of the \emph{inverse} corrections has a direct effect in the precision with which the initialization constants are computed from a given set of initial conditions. Therefore, not computing the inverse transformation with the same accuracy as that of the secular terms will also contribute to the deterioration of the latter. In practice, the highest accuracy of the inverse transformation can be limited to the computation of the periodic corrections of the initial semimajor axis (or its partner canonical variable, the Delaunay action) because it affects directly the secular mean motion $n_\ell$, while the other elements only affect the frequencies of the secular node and perigee, which are $\mathcal{O}(J_2)$ when compared to $n_\ell$.
\par

\subsection{Accuracy of the periodic corrections}

\begin{figure}[htb]
\centering
\includegraphics[scale=0.75]{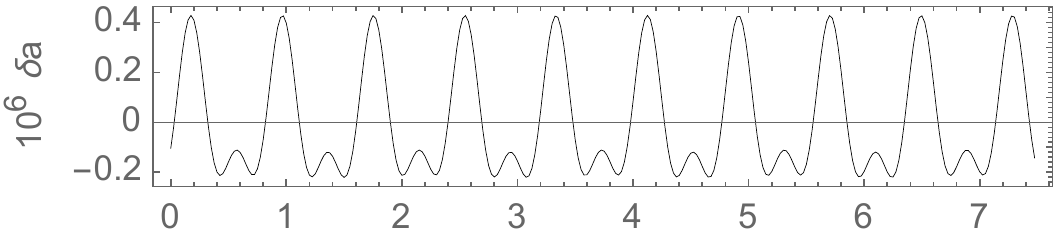}
\includegraphics[scale=0.75]{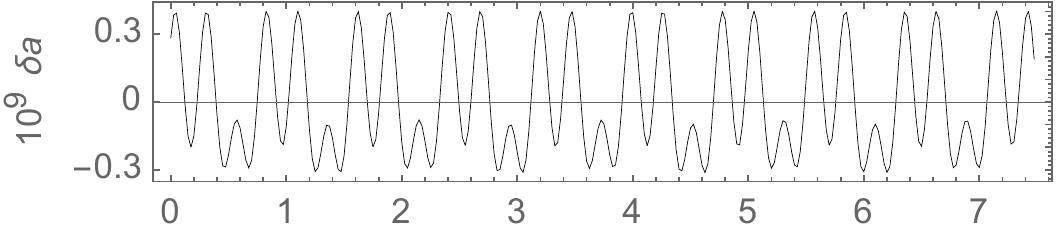}
\includegraphics[scale=0.75]{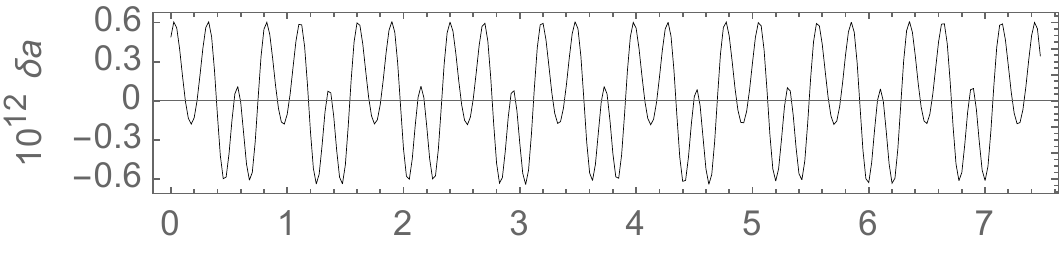}
\includegraphics[scale=0.75]{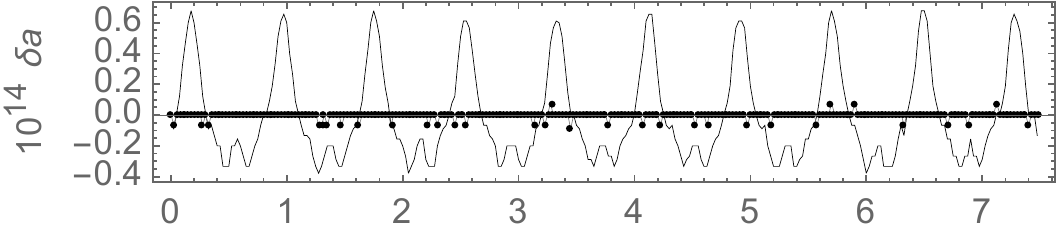}
\caption{Relative errors of the semimajor axis of the PRISMA-type orbit for, from top to bottom, the 1st, 2nd, 3rd, and 4th order truncations of the perturbation solution. The latter shows the relative errors of the 5th order truncation superimposed. Abscissas are hours.}
\label{f:PRISMAdar}
\end{figure}

First of all, we check the accuracy of the periodic corrections for increasing orders of the perturbation solution. If the periodic corrections were exact, transforming different states of the same orbit provided by the reference solution will result in the same, constant secular values of the momenta, and in an exactly linear growing of the secular angles. On the contrary, the secular values obtained from the different states of the true orbit oscillate periodically due to the truncation order of the solution. This is illustrated in Fig.~\ref{f:PRISMAdar} for the semimajor axis of the PRISMA test orbit, where the relative errors with respect to a constant reference value, which has been obtained as the arithmetic mean of the computed secular values, are shown for different truncations of the perturbation solution. We recall that the expected errors of a perturbation solution are of the order of the neglected terms, as follows from Eq.~(\ref{mainPp}). Therefore, for a truncation to the order $m$, we would expect errors of the order $\varepsilon^{m+1}$.
\par

We found that the relative errors of the secular semimajor axis of the first order truncation (top plot of Fig.~\ref{f:PRISMAdar}) are of the order of $J_2^2$, as it should be the case, which amounts to 3 meters in absolute value. The relative errors of the second order truncation  (second to top plot) are of the order of $J_2^3$, or less than 3 millimeter in absolute value. The third and fourth order truncations of the solution yield relative errors of the order of $J_2^4$ and $J_2^5$, respectively (second from bottom and bottom plots, respectively), or absolute errors of micrometers and hundredths of micrometers, respectively. The latter are very close to the 15 exact digits reachable working in double precision. Still, additional improvements are obtained when the truncation of the periodic corrections is extended to the fifth order of $J_2$, now effectively reaching the numerical precision, as shown in the bottom plot of Fig.~\ref{f:PRISMAdar}, in which the relative errors of the fifth order truncation (black dots) are superimposed to the previous case.
\par

\begin{figure}[htb]
\centering
\includegraphics[scale=0.75]{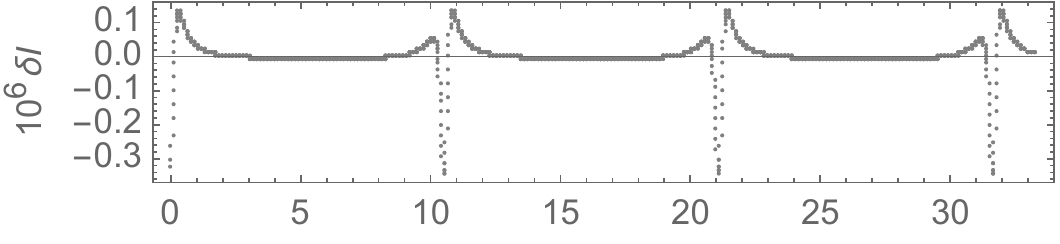}
\includegraphics[scale=0.75]{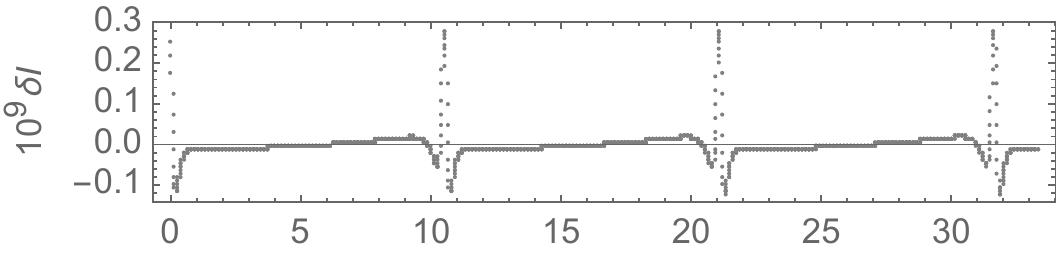}
\includegraphics[scale=0.75]{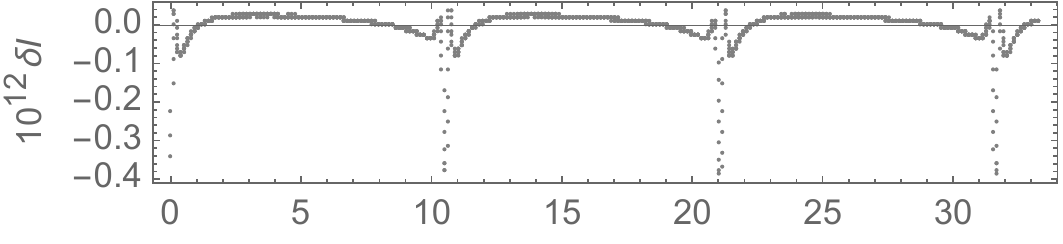}
\includegraphics[scale=0.75]{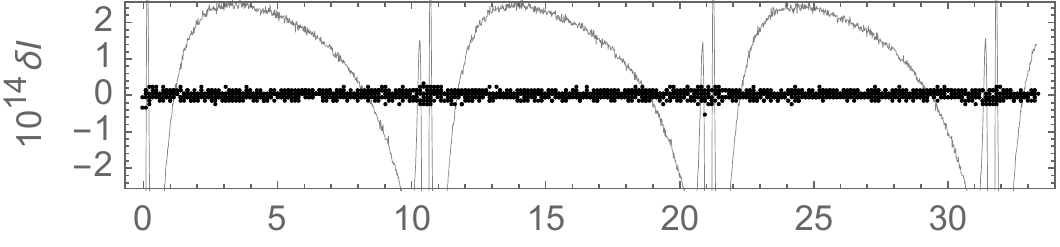}
\caption{Relative errors of the inclination of the GTO orbit for, from top to bottom, the 1st, 2nd, 3rd, and 4th order truncations of the perturbation solution. The latter is superimposed to the relative errors of the 3rd order truncation. Abscissas are hours.}
\label{f:GTOdir}
\end{figure}

The behavior is analogous in the case of the GTO orbit, yet now the errors notable peak at perigee passages. This is illustrated in Fig.~\ref{f:GTOdir}, where, from top to bottom, we present now the relative errors of the inclination $I=\arccos(H/G')$, rather than the total angular momentum $G'$, for the 1st, 2nd, \dots, 4th order truncations of the periodic corrections, respectively. The relative errors of the 4th order truncation in the bottom plot of Fig.~\ref{f:GTOdir} (black dots) are superimposed to the third order ones (gray line) for reference. Corresponding absolute errors are of the order of tens of milliarc seconds for the 1st order, hundredths of mas for the second, hundredths of microarc seconds for the third, and below the level of 1 thousandth of $\mu\mathrm{as}$ for the fourth order. Due to the larger semimajor axis, the GTO orbit is less perturbed in general, and we found that the numerical precision is achieved at the fourth order of the perturbation approach.
\par

\begin{figure}[htb]
\centering
\includegraphics[scale=0.75]{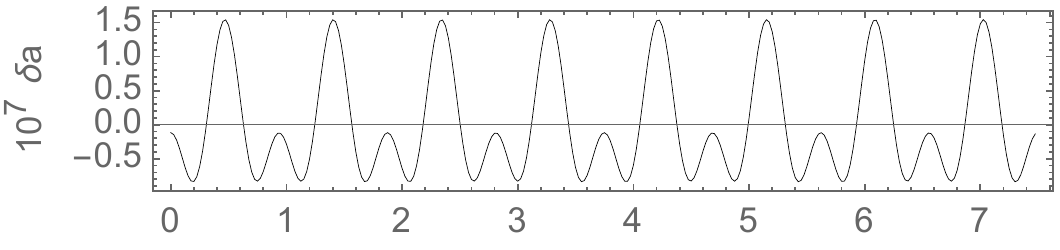}
\includegraphics[scale=0.75]{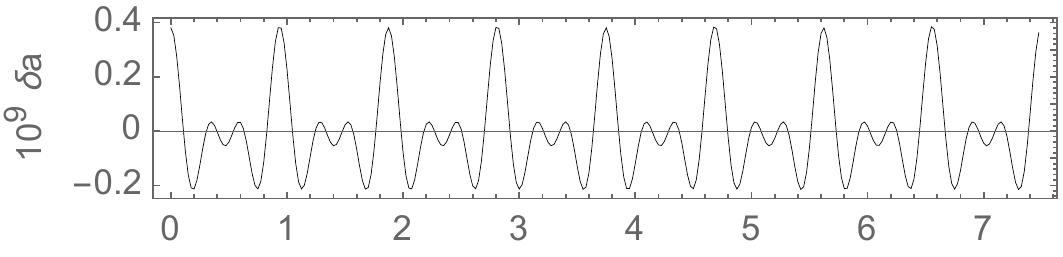}
\includegraphics[scale=0.75]{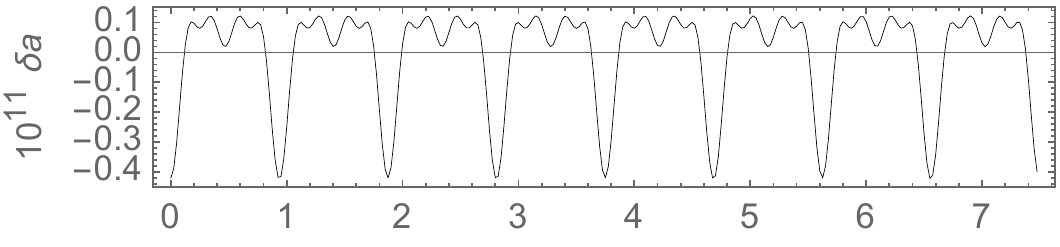}
\includegraphics[scale=0.75]{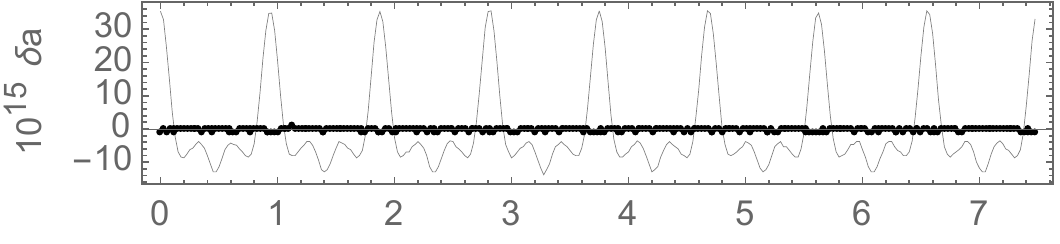}
\caption{Relative errors of the semimajor axis of the TOPEX-type orbit for different truncations of the perturbation solution. Abscissas are hours.}
\label{f:TOPEXdar}
\end{figure}

In the case of the TOPEX orbit the inclination of $\sim66^\circ$ makes the coefficient $5s^2-4$ to become less than 2 tenths. However, this small divisor is not harmful at all in a first order truncation because, as follows from Eq.~(\ref{C1g}), it is multiplied by the square of the eccentricity, which is really small for the TOPEX orbit (about one thousandth, on average). Hence, the relative errors in the semimajor axis introduced by the periodic corrections are slightly better than those previously found for PRISMA, as shown in the top plot of Fig.~\ref{f:TOPEXdar}. The improvements should be a consequence of the larger semimajor axis of the TOPEX orbit when compared with that of PRISMA.
\par

The effects of the offending divisors are more apparent when using higher order truncations of the periodic corrections. Indeed, the coefficients of $\sin(f+2\omega)$ and $\sin(3f+4\omega)$ in Eq.~(\ref{W2g}), as well as the coefficient of $\sin{f}$ in Eq.~(\ref{W2l}), are multiplied by the factor $e/(5s^2-4)^2$, which, while being bigger than before, it is neither of worry for TOPEX yet, as demonstrated in the second to the top plot in Fig.~\ref{f:TOPEXdar}, where it can be checked that the relative errors in the semimajor axis remain of $\mathcal{O}(J_2^3)$. On the contrary, at the third order we find the small divisor alone in coefficients of $\sin2m(f+\omega)$, $m=1,2,3$, in Eq.~(\ref{W3g}), in this way lessening to some extent the corrector effect of the terms factored by $J_2^3$. This slight deterioration of convergence is noted in the second-from-bottom plot of Fig.~\ref{f:TOPEXdar}, where we see that the relative errors in the semimajor axis are now about 10 times bigger than corresponding errors of PRISMA. Bigger offending coefficients are found in higher orders of the periodic corrections, which make that the relative errors of the fourth order truncation are clearly larger than corresponding ones of PRISMA. These relative errors are depicted in the bottom plot of Fig.~\ref{f:TOPEXdar}  (gray line) superimposed with the errors of the fifth order truncation (black dots), showing that, at the end, the numerical precision is reachable also in this case with the fifth order truncation of the solution. The absolute error of the different truncations are now of the order of one meter, a few tenths of cm, sevral hundredths of mm, tenths of $\mu\mathrm{m}$, and thousandths of $\mu\mathrm{m}$, respectively.
\par

\subsection{Accuracy of the analytical solution}

The reference orbits of the three test cases are now compared with the ones provided by the analytical solution. That is, starting from initial conditions in Table \ref{t:testcases}, we first apply the inverse periodic corrections to initialize the constants of the secular solution. Then, the secular terms are evaluated in Eqs.~(\ref{Fsecular})--(\ref{Hsecular}) at the same times $t_i$ as those in which the true solution has been obtained from the numerical integration. It follows the application of the direct periodic corrections to each secular term to get the ephemeris predicted by the analytical solution at the times $t_i$. The accuracy of the different truncations of the analytical solution is assessed by computing the root-sum-square (RSS) of the difference between the position and velocity predicted by the analytical solution and those of the reference orbit and the times $t_i$.
\par

\begin{figure}[htbp]
\centering
\includegraphics[scale=0.75]{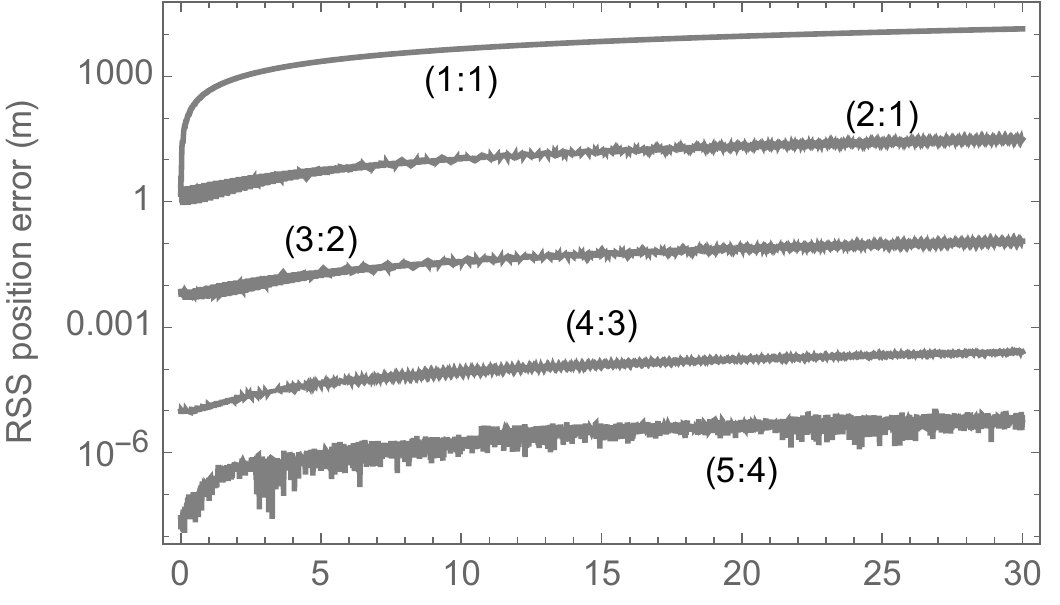}
\caption{RSS errors of different (S:P) truncations of the secular (S) and periodic terms (P) of the analytical solution in the PRISMA test case. Abscissas are days.}
\label{f:PRISMARSS}
\end{figure}

The results of the PRISMA test case are shown in Fig.~\ref{f:PRISMARSS}, in which ordinates are displayed in a logarithmic scale to ease comparisons. The letter $S$ in the labels ($S$:$P$) stand for the truncation order of both the inverse corrections and the secular terms, whereas $P$ indicates the truncation order of the direct periodic corrections. As we already pointed out, the latter do not need to be computed to the same accuracy as the inverse transformation. We observe in Fig.~\ref{f:PRISMARSS} that a simple first order truncation of the analytical solutions ---curve labeled (1:1)--- starts with a RSS error of approximately 1 meter in position, but, due to the errors introduced by the first order truncation of the secular terms, the RSS errors reach more than 10 km after one month. The solution is notably improved when taking the second order terms of the inverse corrections and the secular terms into account. This is shown with the curve labeled (2:1), that only reaches about 30 m at the end of the one-month propagation. The improvements are obtained with only a slight increase of the computational burden, because the inverse corrections and the secular frequencies are evaluated only once. The (3:2) propagation starts from a much smaller RSS error, of less than 1 cm, that only grows to about 10 cm by the end of day 30. Errors fall clearly below the mm level in the case of the (4:3) truncation, and to just a few $\mu\mathrm{m}$ in the (5:4) case. No further improvement of the RSS errors is observed if the direct periodic corrections are taken also up to the fifth order, yet a slight improvement is achieved in that case in the preservation of the energy integral, that reaches in this last case the numerical precision.
\par

\begin{figure}[htbp]
\centering
\includegraphics[scale=0.75]{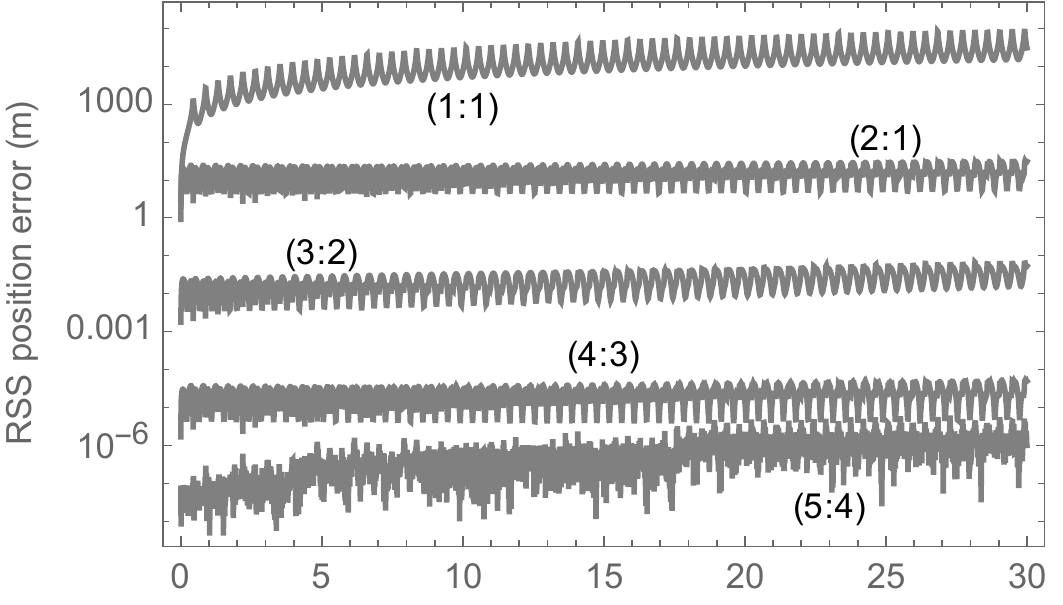}
\caption{RSS errors of different (S:P) truncations of the analytical solution in the GTO test case. Abscissas are days.}
\label{f:GTORSS}
\end{figure}

The behavior of the analytical solution is analogous in the case of the high-eccentricity GTO orbit in what respects to the secular terms, yet the errors of the periodic corrections are now of larger amplitude, as clearly observed in Fig.~\ref{f:GTORSS}. Now, the periodic oscillations of the (1:1) truncation may allow the RSS errors to grow to many tens of km. A secular trend in the RSS errors of about half meter per day of the (2:1) solution remains almost hidden along the 30 days propagation under periodic oscillations of about 30 m. A similar behavior is observed in the (3:2) propagation, where the amplitude of the oscillations is now reduced to the cm order. This amplitude is further reduced to hundredths of mm with the (4:3) truncation. Finally the (5:4) solution improves slightly the propagation, and the RSS errors remain of just a few $\mu\mathrm{m}$ along the whole propagation interval, thus showing the same quality as in the previous test case.
\par

\begin{figure}[htbp]
\centering
\includegraphics[scale=0.75]{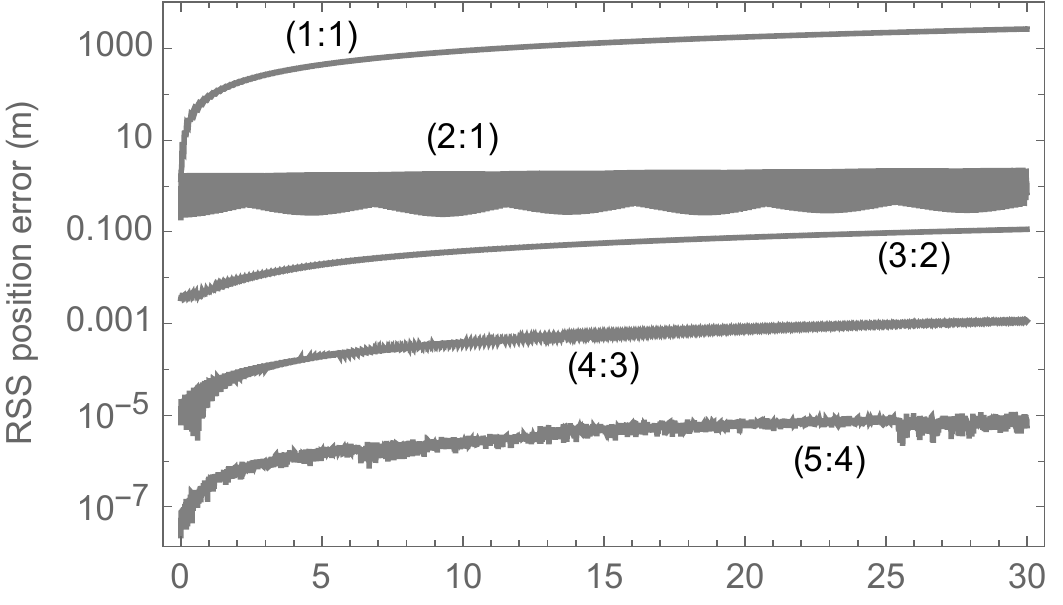}
\caption{RSS errors of different (S:P) truncations of the analytical solution in the TOPEX test case. Abscissas are days.}
\label{f:TOPEXRSS}
\end{figure}

Analogous tests carried out for TOPEX are presented in Fig.~\ref{f:TOPEXRSS}. The results are similar to those previously presented in Fig.~\ref{f:PRISMARSS} for PRISMA, and the apparent discrepancies when using the (2:1) solution are only due to the logarithmic scale, which encompasses a different range in each figure. In fact, the periodic oscillations of the RSS errors are of the same amplitude in both cases ($\sim2$ m), yet the secular trend grows at a low rate of less than 2 cm per day in the case of TOPEX while it does almost two orders of magnitude faster for PRISMA ($\sim1$ m/day). This fact should be attributed to the less perturbed orbit of TOPEX due to its higher altitude. The (3:2) solution provides quite similar RSS errors in both cases, TOPEX and PRISMA, whereas the (4:3) solution performs worse for TOPEX, whose RSS errors grow now at a secular rate about 5 times faster than in the case of PRISMA. These behavior is in agreement with the slight deterioration of the solution previously observed when testing the periodic corrections due to the proximity to the critical inclination. Things are balanced with the (5:4) order solution ---also in agreement with the previously observed behavior due to the increased precision of the analytical solution--- for which the RSS errors reach only the micrometer level at the end of the 30 day propagation.
\par

\begin{figure}[htbp]
\centering
\includegraphics[scale=0.75]{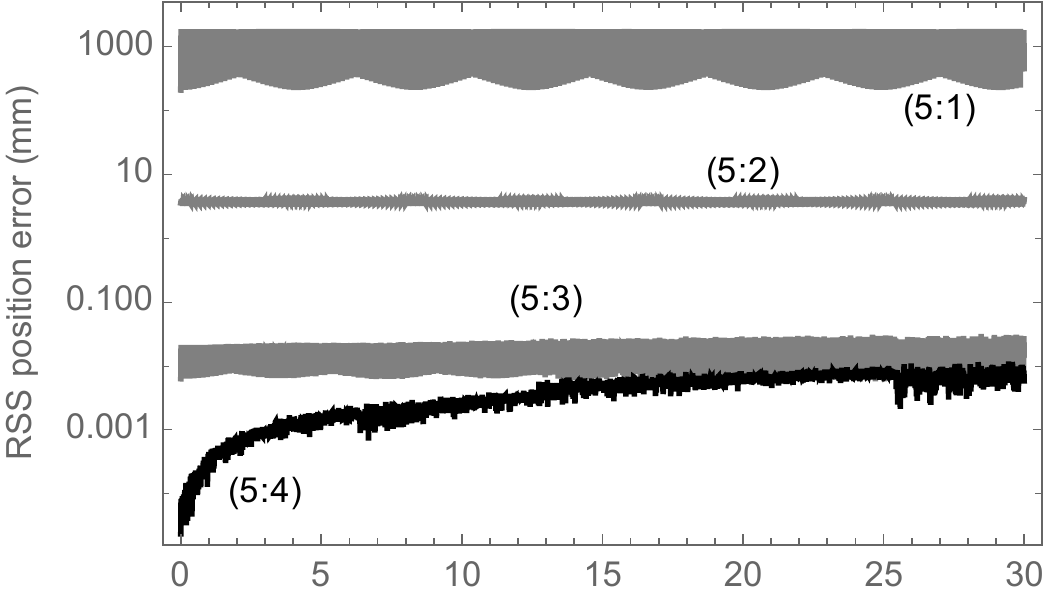}
\caption{TOPEX: RSS errors of different (5:P) truncations. Abscissas are days.}
\label{f:TOPEXRSSb}
\end{figure}

As it is confirmed by the examples carried out, at the end of a long-term propagation the secular trend of the errors prevails over their periodic oscillations. Therefore, in practice, perturbation solutions in which the periodic terms are recovered to a lower order than the order at which the secular frequencies are truncated make full sense. The computational burden of these kinds of solutions is notably alleviated, and hence are definitely much more practicable for ephemeris computation. This is further illustrated for TOPEX in Fig.~\ref{f:TOPEXRSSb}, where it is shown that the (5:3) truncation of the analytical solution might replace the more accurate (5:4) truncation in a long-term propagation with substantial reductions in computing time and minimum degradation of accuracy, which, besides, remains much more uniform along the whole propagation interval. Note that the scale of the ordinates axis has been changed from meters in Fig.~\ref{f:TOPEXRSS} to mm in Fig.~\ref{f:TOPEXRSSb}.
\par

\section{Conclusions}

The main problem of artificial satellite theory is a two degrees of freedom model that captures the bulk of the gravitational effects undergone by low Earth orbits. Still, it lacks of enough integrals to obtain a closed form solution. On the other hand, perturbation methods disclose as specially successful in providing analytical solutions that can replace the true dynamics within a high degree of approximation. Using them, it is shown that there are wide regions of phase space in which the main problem behaves as if it were separable.
\par

In particular, the main problem Hamiltonian has been completely reduced by reverse normalization. That is, the motion of the orbital plane is first decoupled from the motion on that plane, and then the remaining short-period terms are removed by averaging. The normalization is achieved without need of resorting to Hamiltonian simplification procedures, in this way radically departing from a well stablished tradition in artificial satellite theory. The computed analytical solution is exact to all practical purposes in the sense that, working in double-precision floating point aritmethic, there are no substantial differences with respect to exact numerical integrations (computed with extended precision) for reasonably long time intervals. In agreement with the particular value of the Earth's oblateness coefficient, this precision is achieved only when the perturbation approach is extended up to, at least, the fifth order of the small parameter.
\par

\begin{acknowledgements}
Support by the Spanish State Research Agency and the European Regional Development Fund under Projects ESP2016 -76585-R and ESP2017-87271-P (MINECO/ AEI/ERDF, EU) is recognized.

\paragraph*{\small\bf Conflict of Interest}: The author declares that he has no conflict of interest.

\end{acknowledgements}

\end{document}